\documentclass[12pt]{article}

\usepackage{algorithm}
\usepackage{algorithmic}
\usepackage{amsmath,amscd}
\usepackage{amsfonts}
\usepackage{amssymb}
\usepackage{times}
\usepackage{ulem}
\usepackage{comment}
\usepackage{graphicx}
\usepackage{float}
\usepackage{subcaption}
\usepackage{color}
\usepackage{url}
\graphicspath{{figures/}}

\usepackage{geometry}
\newcommand{\blue}[1]{{\color{blue}{#1}}}

\title{Entrainment dynamics of forced hierarchical circadian systems revealed by 2-dimensional maps}

\author
{Guangyuan Liao, Casey Diekman, Amitabha Bose\\
	\\
	\normalsize{Department of Mathematical Sciences, New Jersey Institute of Technology, Newark, NJ 07102}\\
}

\date{\today}

\begin{document}
	
	
	\baselineskip24pt
	
	
	\maketitle
	
	
	\begin{abstract}
		The ability of a circadian system to entrain to the 24-hour light-dark cycle is one of its most important properties.  A new tool, called the entrainment map, was recently introduced to study this process for a single oscillator. Here we generalize the map to study the effects of light-dark forcing in a hierarchical system consisting of a central circadian oscillator that drives a peripheral circadian oscillator. We develop techniques to reduced the higher dimensional phase space of the coupled system to derive a generalized 2-D entrainment map. Determining the nature of various fixed points, together with an understanding of their stable and unstable manifolds,  leads to conditions for existence and stability of periodic orbits of the circadian system. We use  the map to investigate how various properties of solutions depend on parameters and initial conditions including the time to and direction of entrainment. We show that the concepts of phase advance and phase delay need to be carefully assessed when considering hierarchical systems.
	\end{abstract}
	
	
	
	
	\section{Introduction}
	\label{sec:intro}
	
	Circadian rhythms refer to a variety of oscillatory processes that occur  over a roughly 24-hour time period. Circadian oscillations are found in a variety of animal and plant species \cite{Ben-ShlomoKyriacou02}. Within humans  a common  example involves our core body temperature which shows a local minimum typically in the early morning hours ($\sim$ 4:00 AM) and a local maximum roughly twelve hours later \cite{krauchi2002circadian}. Similarly, concentrations of certain hormone levels within our bodies oscillate over the course of a day \cite{czeisler1999circadian}. In the absence of any explicit forcing from naturally occurring light-dark cycles, circadian oscillators possess endogenous periods of roughly 24 hours. Their ability to also entrain to  24-hour periodic cycles of light and dark is one of their most important properties.
	
	The entrainment of circadian oscillators has been mathematically analyzed using a variety of techniques. Often this involves describing the circadian oscillator with a reduced phase description such as that given by a Kuramoto oscillator \cite{bordyugov2015tuning,Kuramoto84}. The problem then reduces to studying periodically forced Kuramoto systems. Other approaches include deriving model equations that retain more of their connection to the underlying biological process \cite{tyson1999simple,novak2008design}. Recently Diekman and Bose \cite{diekman2016entrainment} introduced a novel tool called the entrainment map to determine whether  a circadian oscillator can entrain to the 24-hour light-dark cycle, and if so, at what phase.  The derived map is equivalent to a 1-D Poincar\'e map that tracks the phase of light onset of the light-dark forcing on a cycle-by-cycle basis. In principle, the  dimension of the underlying circadian oscillator model is not relevant. Diekman and Bose derived entrainment maps for the 2-D Novak-Tyson model \cite{tyson1999simple,novak2008design}, the 3-D Gonze model \cite{GonzeHerzel05} and the 180-D Kim-Forger model \cite{KimForger12}. In general,  the map can be used to estimate both entrainment times and whether entrainment occurs through phase advance or delay with respect to the daily onset of lights.
	
	There are several scenarios in which circadian oscillators do not directly receive light-dark forcing \cite{GonzeHerzel05,gu2017entrainment, LeiseSiegelmann06}. Instead they are part of hierarchical systems in which, as ``peripheral'' oscillators,  they are periodically forced by other ``central''  circadian oscillators that do directly receive light input. Cells within major organs in our bodies fall into this category. Several natural questions arise about the entrainment process of these peripheral oscillators. For example, do they entrain through phase advance or phase delay as central oscillators do? To what extent is their entrainment time dependent on the entrainment process of the central oscillator from which they receive forcing? To study such questions, here we generalize the entrainment map to a 2-D map where we track from the perspective of the peripheral oscillator both the phase of the central oscillator as well as the phase of light onset.
	
	In this paper, we first consider the situation in which a single central oscillator receives light-dark input. In turn, this central oscillator sends input to a single peripheral oscillator. To focus on the mathematical aspects of the derivation and analysis of the 2-D entrainment map, we will utilize the planar Novak-Tyson model \cite{novak2008design} for both the central and peripheral oscillators. The phase space for this problem is 5-D, two for each of the oscillators and a fifth that accounts for the light dark forcing.  We will define a Poincar\'e section transversal to the flow allowing us to derive a 2-D map that determines the phase of light and the phase of the central oscillator at each cycle when the peripheral oscillator lies on the Poincar\'e section. We analyze the map by extending techniques first introduced in Akcay et al \cite{akcay2014effects,akcay2018phase}. We will show that for a range of parameter values, the map possesses four fixed points:  one asymptotically stable and three unstable fixed points, two of which are saddle points. All of these fixed points are related to actual periodic orbits of the flow. By numerically calculating entrainment times (defined precisely later in the text), we are able to uncover how the stable and unstable manifolds of the saddle points organize the iterates of the map, determine the direction of entrainment and give rise to a rich set of dynamics The findings of the map are then validated by comparing them to direct simulations of the model equations. We also extend the analysis to the case of a semi-hierarchical system that consists of a second central oscillator that receives less light input than the first central oscillator.
	
	Analysis of the map reveals several important insights into the entrainment and reentrainment process. First, bounds on important parameters, such as the intensity of light input and the strength of the coupling from the central oscillator that lead to entrainment, are easily identified. We are able to determine which kinds of perturbations lead to faster or slower reentrainment, e.g. whether perturbations that desynchronize only the peripheral oscillator but not the central one lead to quick reconvergence. Interestingly, we find that the straightforward notion of convergence via phase advance or phase delay needs to be generalized. Indeed, the peripheral oscillator can converge by a combination of phase advance and delay while central oscillators typically converge by either phase advancing or delaying. The direction of entrainment is related to the concept of orthodromic versus antidromic reentrainment, which in circadian systems,  can lead to ``convergence by partition" \cite{LeiseSiegelmann06}, as will be elaborated upon in the Discussion.

	\section{Models and Methods}
	\label{sec:Models}
	
	Our model is based on the Novak-Tyson (NT) model \cite{tyson1999simple} for the molecular circadian clock in the fruit fly {\it Drosophila}. The NT model can be written in the following form:
	\begin{align}
	\begin{split}
	\frac{1}{\phi}\frac{dP}{dt} &= M - k_{f}h(P) - k_{D}P  - k_{L}f(t)P\\
	\frac{1}{\phi}\frac{dM}{dt} &= \epsilon\left(g(P) - M\right)
	\end{split}
	\end{align}
	where $g(P) = \frac{1}{1+P^4}$, and $h(P) = \frac{P}{0.1+P+2P^2}$ . The $M$ variable represents mRNA concentration, and $P$ variable represents the protein concentration. The parameter $\epsilon$ is small, which separates $P$ and $M$ into fast and slow variables. The parameter $\phi$ will directly affect the period of the solutions of this system; smaller values imply longer endogenous periods. The function $f(t)$ describes the light-dark (LD) forcing, which is defined by a 24-hour periodic step function, $f(t)= 1$ when lights are on and $f(t)=0$ when lights are off. We consider for convenience  a 12:12 photoperiod though there is no difficulty in extending to other cases. In {\it Drosophila}, there is protein degradation during darkness, and light increases the degradation. So $k_D$ represents the degradation rate during darkness, and $k_L$ represents the additional degradation rate which is caused by light. The parameter $k_f$ is a combination of two variables in the original Novak\b;ue{-}Tyson paper \cite{tyson1999simple}.
	In \cite{diekman2016entrainment}, the entrainment of a single NT oscillator to a 24-hour period LD forcing was studied. The ensuing solution was denoted as an LD-entrained solution.
	
	\subsection{Coupled Novak-Tyson Model}
	The coupled Novak-Tyson (CNT) model is given by the following equations:
	\begin{align}
	\begin{split}
	\frac{1}{\phi_1}\frac{dP_1}{dt} &= M_1-k_{f}h(P_1)- k_{D}P_1 - k_{L_1}f(t)P_1 \\
	\frac{1}{\phi_1}\frac{dM_1}{dt} &= \epsilon [g(P_1) -M_1 ]	\\
	\frac{1}{\phi_2}\frac{dP_2}{dt} &= M_2- k_{f}h(P_2)- k_{D}P_2 - k_{L_2}f(t)P_2 \\
	\frac{1}{\phi_2}\frac{dM_2}{dt} &= \epsilon [g(P_2) -M_2 + \alpha_1 M_1 g(P_2)]
	\end{split}
	\label{CNT}
	\end{align}
	The parameters and variables have the same meaning as the original NT model. We introduce a coupling term $\alpha_1 M_1g(P_2)$, from oscillator 1 ($O_1$) to oscillator 2 ($O_2$). The parameter $\alpha_1$ is a non-negative real number which denotes the coupling strength. We put the coupling factor into the second equation of $O_2$ based on Roberts et al. \cite{roberts2016functional}, who suggest that coupling occurs between the mRNA production rates.
	
	To find the entrained solutions and understand the geometry of the CNT system in the presence of the LD cycle, the nullclines of each oscillator play an important role. The nullclines are the set of points where the the right hand sides of (\ref{CNT}) equal zero. We define the $P$-nullcline and $M$-nullcline for both oscillators as follows:
	\begin{align}
	\begin{split}
	N_{P_D}: M&= k_f h(P) + k_D P  \\
	N_{P_L}: M&= k_f h(P) + (k_D + k_L)P\\       
	N_{M_1}:M_1&= g(P_1)\\
	N_{M_2}:M_2&= g(P_2) + \alpha_1 M_1 g(P_2)
	\end{split}
	\end{align}
	where the $P$-nullcline are cubic shaped curves for both oscillators. Since the LD forcing is a Heaviside function that switches between 0 and 1, $N_{P_D}$ represents the $P$-nullcline in constant dark condition and $N_{P_L}$ is for constant light condition.  $N_{M_1}$ is a sigmoidal like curve,  as is $N_{M_2}$. However the latter can vary depending on the value of $M_1$. When $O_1$ is entrained, along its limit cycle, the $M_1$ value is bounded between $min|M_1(t)|$ and $max|M_1(t)|$. Thus $N_{M_2}$ can oscillate between $g(P_2) + \alpha_1 min|M_1(t)| g(P_2)$ and $g(P_2) + \alpha_1 max|M_1(t)| g(P_2)$. We assume that any intersection between $N_{P}$ and $N_{M}$ occurs on the middle branch of the corresponding cubic nullclines. This will guarantee that any ensuing fixed points of the CNT system are unstable and will allow oscillations to exist. The nullclines are shown in Fig. \ref{fig:O1_limitcycle} and \ref{fig:O2_limitcycle}.
	
	\begin{figure}[!htbp]
		\centering
		\begin{subfigure}{.45\linewidth}
			\centering
			\includegraphics[width=.7\linewidth]{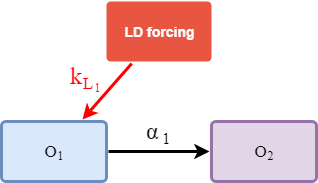}
			\caption{}
			\label{fig:case1}
		\end{subfigure}
		\begin{subfigure}{.45\linewidth}
			\centering
			\includegraphics[width=.7\linewidth]{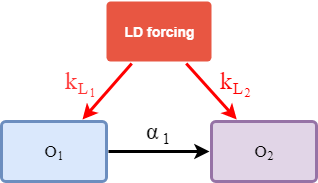}
			\caption{}
			\label{fig:case_nofeedback}
		\end{subfigure}
		\caption{(a): Model with strict hierarchical coupling. (b): Semi-hierarchical model when both oscillators receive light input, but the light into $O_2$ is much weaker than the light into $O_1$.}
	\end{figure}
	
	We mainly study the case with strict hierarchical coupling, which is shown in Fig. \ref{fig:case1}. In this case, the LD forcing is applied only on $O_1$, which then has feedforward coupling onto $O_2$. We fix the value of parameter $k_{L_2}=0$. Figure \ref{fig:case_nofeedback} shows the semi-hierarchical CNT model when both oscillators receive light forcing, but the effect of light into $O_2$ is taken to be less intense than that into $O_1$, in other words, $k_{L_2}<k_{L_1}$.
	
	\subsection{The entrainment map}
	
	When attempting to determine the existence of periodic solutions using Poincar\'e maps, one has to decide where in phase space to place the section. Often in circadian models, the Poincar\'e section is placed on the 24-hour light-dark forcing, leading to a stroboscopic map that determines the state of the system every 24 hours. In \cite{diekman2016entrainment}, Diekman and Bose instead placed the section in the phase space of the circadian oscillator and backed out the phase of light when the oscillator was at the section. Here, we have followed that approach when building the 2-D map.  The  Poincar\'e section was chosen at a location in the flow that $O_2$ crossed. The phase of $O_1$ with respect to a reference point on its own limit cycle, $x$, and of lights $y$ was then determined to derive the 2-D map. In this section, We first introduce the original 1-D map, and then generalize it to our 2-D map.
	
	The entrainment map $\Pi(y)$ for the original NT model was introduced as a 1-D map in \cite{diekman2016entrainment}. To define $\Pi(y)$, Diekman and Bose take a Poincar\'e section $\mathcal{P}$ as a 1-D line segment which intersects the LD-entrained solution of a single periodically forced NT oscillator. The section is chosen along the left branch of $N_p$ at $M=0.45, |P-0.852|<\delta$, where $\delta$ is a small control variable such that the section is fixed in a neighborhood of the $P$-nullcline. This location is a natural choice since all trajectories of the NT oscillator pass through the section.  A phase variable $y$ is defined to be the amount of time that has passed since the beginning of the most recent LD cycle. When the trajectory first returns to $\mathcal{P}$, the map $\Pi(y)$ is defined to be the amount of time that has passed since the onset of the most recent LD cycle, which is the new phase of the light forcing. The domain and range of $\Pi(y)$ are both (0,24]. The domain is actually homeomorphic to the unit circle $\mathbb{S}^1$,
	so $y=0$ and $y=24$ are equivalent. The map is written as $y_{n+1}=\Pi(y_n)$, where:
	
	\begin{equation}
	\Pi(y_n)=(\rho(y_n)+y_n)\ mod\ 24.
	\label{eq:1D map}
	\end{equation}
	$\rho(y)$ is a return time map that measures the time a trajectory starting on $\mathcal{P}$ takes to return to $\mathcal{P}$. It is continuous and periodic at its endpoints $\rho(0^+)=\rho(24^-)$. If $\rho(y)<24-y$, then $\Pi(x)=\rho(y)+y$, because the trajectory will return back to $\mathcal{P}$ within the same LD cycle which it started. If $24-y<\rho(y)<48-y$, then $\Pi(y)=\rho(y)+y-24$, because the trajectory will return in the next LD cycle and so on.
	
	If there exists a $y_s$, such that $y_s=\Pi(y_s)$ and $\lvert\Pi'(y_s)\vert<1$, then $y_s$ is a stable fixed point of the map $\Pi(y)$, and it also determines a 1:1 phase locked solution. The phenomenon of 1:1 phase locking in this case occurs when the oscillator has one return to the Poincar\'e section for every one period of the LD forcing. When a stable solution exists, the map $\Pi(y)$ quite accurately calculates the time to approach the stable solution starting from any initial condition of $y$. Numerically we use the concept of entrainment to evaluate the convergence time. Suppose {$y_j$} is a sequence of iterates of the map, then we say the solution is entrained if there exists $m$, such that for all $j\geq m$, $\lvert y_s-y_j\rvert <0.5$. The entrainment time is then $\Sigma_{i=1}^m \rho(y_i)$.
	
	\paragraph{The 1-D $O_1$-entrained map for the CNT system}
	
	The 1-D map for the NT system can not be directly applied to the CNT system, because the second oscillator will have additional free variables to determine, meaning that the entrainment map for the CNT system will be higher dimensional. However, for the hierarchical CNT system, if we assume that $O_1$ is already entrained, then the chain \textit{$LD \Rightarrow O_1 \Rightarrow O_2$} is reduced to \textit{$O_1$}-entrained \textit{$\Rightarrow O_2$}. The system can be rewritten in the following manner:
	\begin{align}
	\begin{split}
	\frac{1}{\phi_2}\frac{dP_2}{dt} &= M_2 - k_{f} h(P_2)-k_{D}P_2 \\
	\frac{1}{\phi_2}\frac{dM_2}{dt} &= \epsilon [g(P_2)-M_2 + \alpha_1 M_1 g(P_2)]
	\end{split}
	\label{Reduced_CNT}
	\end{align}
	In the $O_1$-entrained case, $O_2$ is forced by the coupling from $O_1$, the major difference between the forcing via coupling and via light is that The direct forcing is a heaviside periodic function, in other words, a square wave, but the forcing via coupling is a continuous wave, and it’s not periodic until O1 is entrained. we take a Poincar\'e section that intersects the entrained $O_2$ limit cycle solution at  $\mathcal{P}: P_2=1.72, |M_2-0.1289|<\delta$ such that $P_2'<0$, where $\delta$ is a small control parameter such that all trajectories starting on the section return to the section. See remark for more detail of the choice of Poincar\'e section. After taking the Poincar\'e section $\mathcal{P}$, $P_2$ is fixed, and $M_2$ is bounded by $\delta$, so the only free variable is the phase of light. We define the 1-D $O_1$-entrained map by
	\begin{equation}
	y_{n+1} = \Pi_{O_1}(y_n)=(y_n + \rho(y_n;\gamma(y_n)))\ mod\ 24
	\label{Eq:premap}
	\end{equation}
	where $y\in (0,24]$ is defined to be the phase of the LD forcing, which has the same meaning as the 1-D entrainment map in \cite{diekman2016entrainment}. We define $\gamma(t):=\varphi_t(X_0)$ to be the LD-entrained limit cycle of $O_1$, where $X_0$ is a chosen reference point on $\gamma(t)$. We denote the set of points that lie on the limit cycle of $O_1$  by $\Gamma_{O_1}$. At $X_0$, the lights just turn on for $O_1$. In the $O_1$-entrained case, the location of $O_1$ only depends on $y_n$ and can be denoted by
	$\gamma(y_n)$. Based on the above definition, $\gamma(y_n)$ means a point on the limit cycle of $O_1$ when the light has been turned on for $y_n$ hours. $\rho(y_n)$
	measures the return time when $O_2$ first returns $\mathcal{P}$.
	
	Notice that in the definition of the $O_1$-entrained map, the phase of $O_1$ is determined by $y$ (the phase of the LD forcing), since it is $O_1$-entrained. This makes the $O_1$-entrained map a 1-D map, and most of the properties of the NT model's 1-D map carry over to the $O_1$-entrained map. For example, if there is a point $y_s$, such that $y_s=\Pi_{O_1}(y_s)$ and $|\Pi_{O_1}'(y_s)|<1$, then $y_s$ is a stable fixed point of the $O_1$-entrained map. The fixed points of the map also determine 1:1 phase locked solutions of the coupled system.\\
	
	\paragraph{The general 2-D entrainment map}
	
	In the case of the $O_1$-entrained map, the initial location of $O_1$ when $O_2$ lies on $\mathcal{P}$ is always determined by $y$, the phase of the LD cycle. But in general, the initial location of $O_1$ doesn't always depend on $y$, rather it could lie arbitrarily in its phase space. To limit the possibilities, we restrict \blue{the initial location of} $O_1$ to lie anywhere along its own limit cycle $\Gamma_{O_1}$. This restriction will therefore only introduce one new free variable and motivates us to generalize the map to two dimensions.
	\begin{equation*}
	(x_{n+1},y_{n+1}) = \Pi(x_n,y_n) = (\Pi_1(x_n,y_n),\Pi_2(x_n,y_n))
	\end{equation*}
	We keep the definition of $y_n$ and the location of the Poincar\'e section $\mathcal{P}$ the same as the $O_1$-entrained map. We now introduce a new variable $x$ to determine $O_1$'s position in phase space relative to its own LD-entrained solution. The detailed definition is explained using a phase angle.
	
	\begin{figure}[!htbp]
		\centering
		\begin{subfigure}{.45\linewidth}
			\centering
			\includegraphics[width=\linewidth]{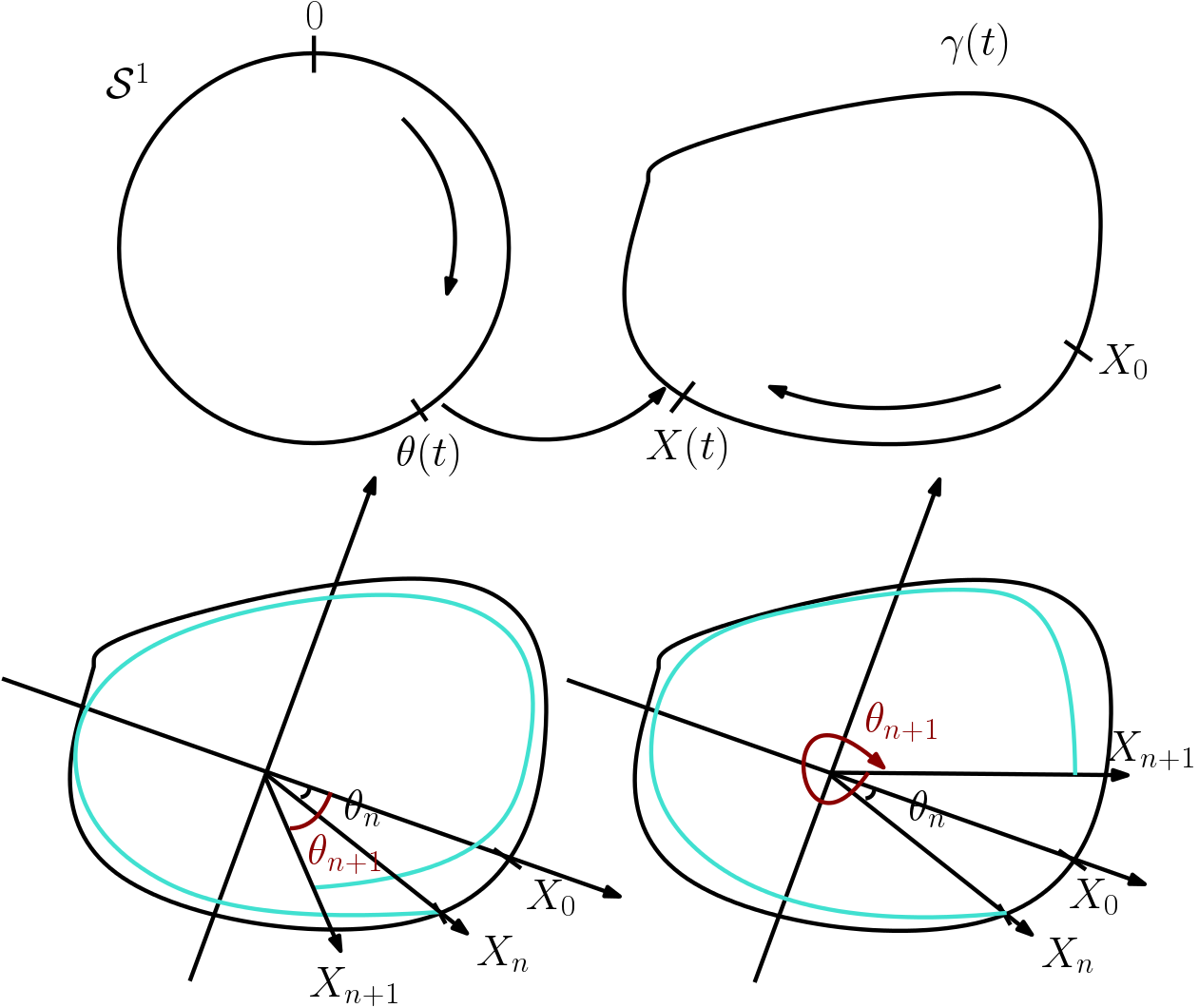}
			\caption{Schematic for $\Pi_1$}
			\label{fig:pi1}
		\end{subfigure}
		\begin{subfigure}{.45\linewidth}
			\centering
			\includegraphics[width=\linewidth]{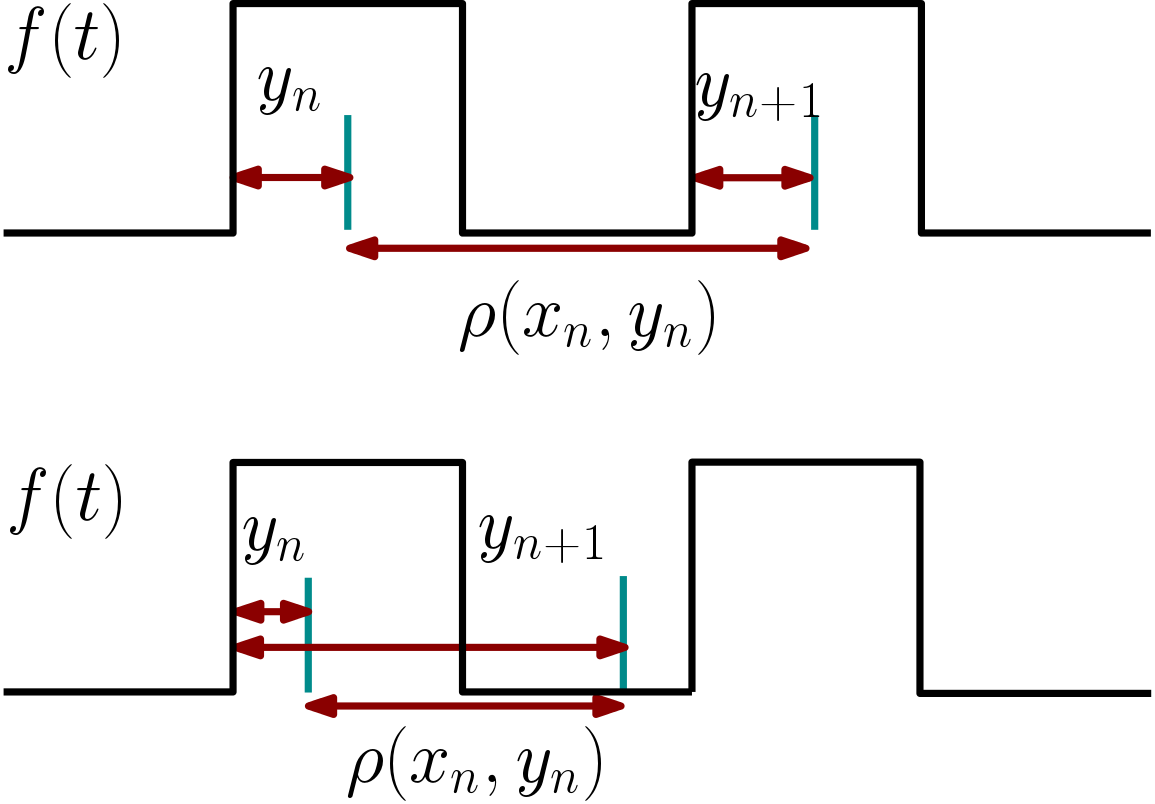}
			\caption{Schematic for $\Pi_2$}
			\label{fig:pi2}
		\end{subfigure}
		\caption{(a): The upper panel shows a schematic of the homeomorphism from the unit circle $\mathbb{S}^1$ to $\Gamma_{O_1}$. The lower panel shows how we construct the map in two different conditions; the left one shows the case when the phase angle $\theta$ associated with the trajectory of $O_1$ rotates through more than 2$\pi$, the right one is where the rotation is less than 2$\pi$. (b): In both panel schematics, the first blue vertical line segment denotes where we chose the initial phase of light. After time $\rho(x_n,y_n)$, the trajectory returns to $\mathcal{P}$, and the new phase of light is $y_{n+1}$. For the upper panel, $y_n+\rho(x_n,y_n)>24$, so $y_{n+1}=y_n+\rho(x_n,y_n)-24$. For the lower panel, $y_n+\rho(x_n,y_n)<24$, so $y_{n+1}=y_n+\rho(x_n,y_n)$. The black square wave in both panel represents the LD forcing.}
		\label{fig:ex_Pi}
	\end{figure}
	\paragraph{Defining $\Pi_1$ using a phase angle}
	According to the $O_1$-entrained map, the  trajectory of $O_1$ always remains on $\Gamma_{O_1}$. However, if $O_1$ is not already entrained, then its trajectory may not lie on $\Gamma_{O_1}$ but will instead approach it asymptotically. Thus we need a new independent variable to determine the position of $O_1$ for this situation. From the $O_1$-entrained case, the position of $O_1$ can always be described as $\gamma(t)$, where $t\in(0,24]$. The idea is to define a new independent phase variable $x$ equivalent to the time variable $t$ that is obtained by projecting the real location of $O_1$ onto its limit cycle $\Gamma_{O_1}$, while keeping the error small.
	
	We define the phase angle in the following steps:
	
	1) Transform the coordinate system appropriately: Shift the origin to the intersection point of the uncoupled $O_1$'s two nullclines $N_P$ and $N_M$. Then connect the origin and the point $X_0$ and expand the line segment as the x-axis of the new coordinate system. The y-axis is determined automatically to be orthogonal to the x-axis, as in Fig. \ref{fig:pi1}.
	
	2) Define $x$ in terms of the phase angle: Consider the phase plane as a complex plane $\mathbb{C}$. Let's call the point $X_0$ as $z_0=r_0e^{i\theta_0}\in \mathbb{C}$, where $\theta_0=0$ after the coordinate system transformation. We can then represent any point on the limit cycle $\gamma(t)$ as a complex number $z=re^{i\theta}$, where we define $\theta\in (0,2\pi]$. Then $x$ is defined to be the phase of $O_1$ when choosing $X_0$ as the reference point. In other words, $z=\gamma(x)=re^{i\theta}$. Notice that $x$ is homeomorphic to the unit circle
	$\mathbb{S}^1$, because $\theta=Arg(\gamma(x))$; see  Fig. \ref{fig:pi1}.
	The domain of $x$ is also $\mathbb{\tilde{S}}^1=(0,24]$.
	
	3) Define the map $\Pi_1$. Suppose we start integrating the system with any initial condition $(x_n,y_n)$  (see Fig. \ref{fig:pi1}, lower panel as an example). After the time $\rho(x_n,y_n)$, $O_2$ returns to the Poincar\'e section, the new location of $O_1$ is now
	\begin{equation*}
	\Psi_{\rho(x_n,y_n)}(\gamma(x_n))=r_{n+1} e^{i\theta_{n+1}}
	\end{equation*}
	where $\Psi_t(X)$ is the flow of $O_1$, and the phase angle is $\theta_{n+1}$. We then find the unique point $\hat{x}$ lying on $\Gamma_{O_1}$ such that the phase angle of $\Psi_{\rho(x_n,y_n)}(\gamma(x_n))$ matches the angle associated with $\gamma(\hat{x})$. That is we choose $\hat{x}$ such that  $Arg(\gamma(\hat{x}))=\theta_{n+1}$. Geometrically, we are simply choosing $\hat{x}$ as the associated value at which  the ray passing through $\Psi_{\rho(x_n,y_n)}(\gamma(x_n))$ intersects $\gamma(t)$.
	We define $x_{n+1}=\hat{x}$. We can then write $\Pi_1$ as the following:
	\begin{equation}
	x_{n+1} = \Pi_1(x_n,y_n) = \{\hat{x} \in [0,24): Arg(\gamma(\hat{x}))=\theta_{n+1}\}.
	\end{equation}
	The numerical algorithm is as following:
		\begin{algorithm}[H]
			\caption{Calculate $\Pi_1$(Initial phase of $O_1$: $x_0$, initial phase of light: $y_0$, $O_1$ limit cycle: $\gamma_1(t)$)}
			\begin{algorithmic}
				\STATE $(P_1,M_1)=\gamma_1(x_0)$
				\STATE $(P_2,M_2)=\mbox{point on Poincar\'e section}$
				\STATE Integration initial = $(P_1,M_1,P_2,M_2,y_0)$
				\STATE Integrate CNT system until $P_2$ returns section
				\STATE Get the ending location $(P_1',M_1')$
				\STATE $\theta = FindAngle(P_1',M_1')$ (calculate the phase angle of the new point)
				\STATE $x_0' = compare(\gamma_1,\theta)$ (Find $\gamma_1(x_0')$, such that $|Arg(\gamma_1(x_0')) - \theta|<\epsilon$
				\RETURN $x_0'$
			\end{algorithmic}
		\end{algorithm}
	
	The definition of $\Pi_2$ is straightforward. We just mimic the construction of the $O_1$-entrained map. The only difference is that the return time function $\rho$ depends on both $x$ and $y$, because $O_1$ is no longer $O_1$-entrained:
	\begin{equation}
	y_{n+1} = \Pi_2(x_n,y_n) = y_n + \rho(x_n,y_n)\ mod\ 24,
	\end{equation}
	where $y \in \mathbb{\tilde{S}}^1=(0,24]$ is defined on a homeomorphism of the unit circle $\mathbb{S}^1$, $y=h(\theta)=\frac{12}{\pi} \times \theta$.
	
	The schematic Fig. \ref{fig:pi2} depicts a way to understand the definition of $\Pi_2$. The first blue vertical line segment signifies the initial phase of $O_2$ when it starts on $\mathcal{P}$ when the light turns on $y_n$ hours. After time $\rho(x_n,y_n)$, the trajectory returns to $\mathcal{P}$, signified by the second blue vertical line segment, with the lights having turned on $y_{n+1}$ hours ago. In the upper panel, $\rho(x_n,y_n)>24-y_n$, therefore the trajectory does not return to $\mathcal{P}$  within the same LD cycle. In the lower panel, $\rho(x_n,y_n)<24-y_n$, therefore the trajectory does return to $\mathcal{P}$ within the same LD cycle.

	\section{Results}
	\label{sec:results}
	
	In this section, we first show simulations demonstrating the entrainment of the strictly-hierarchical CNT model. We then define and analyze a 1-D map in which $O_1$ is assumed to already be entrained. We call this the $O_1$-entrained map. Understanding the 1-D map will facilitate the definition and analysis of the 2-D entrainment map. Finally, we extend the results to the semi-hierarchical case.
	
	\subsection{The entrained solutions of the CNT model}
\begin{figure}[H]
		\centering
		\begin{subfigure}{.45\linewidth}
			\centering
			\includegraphics[width=\linewidth]{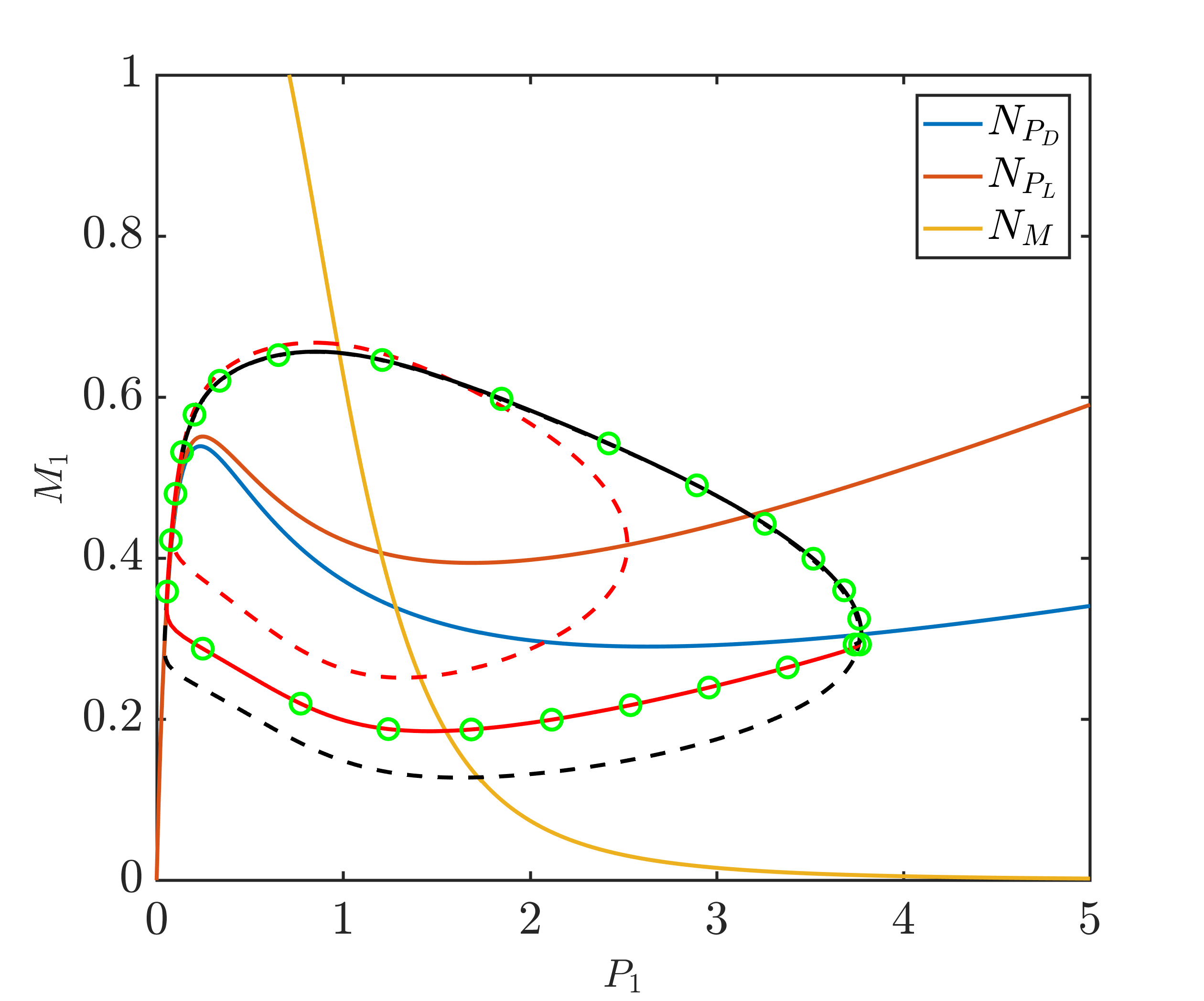}
			\caption{}
			\label{fig:O1_limitcycle}
		\end{subfigure}
		\begin{subfigure}{.45\linewidth}
			\centering
			\includegraphics[width=\linewidth]{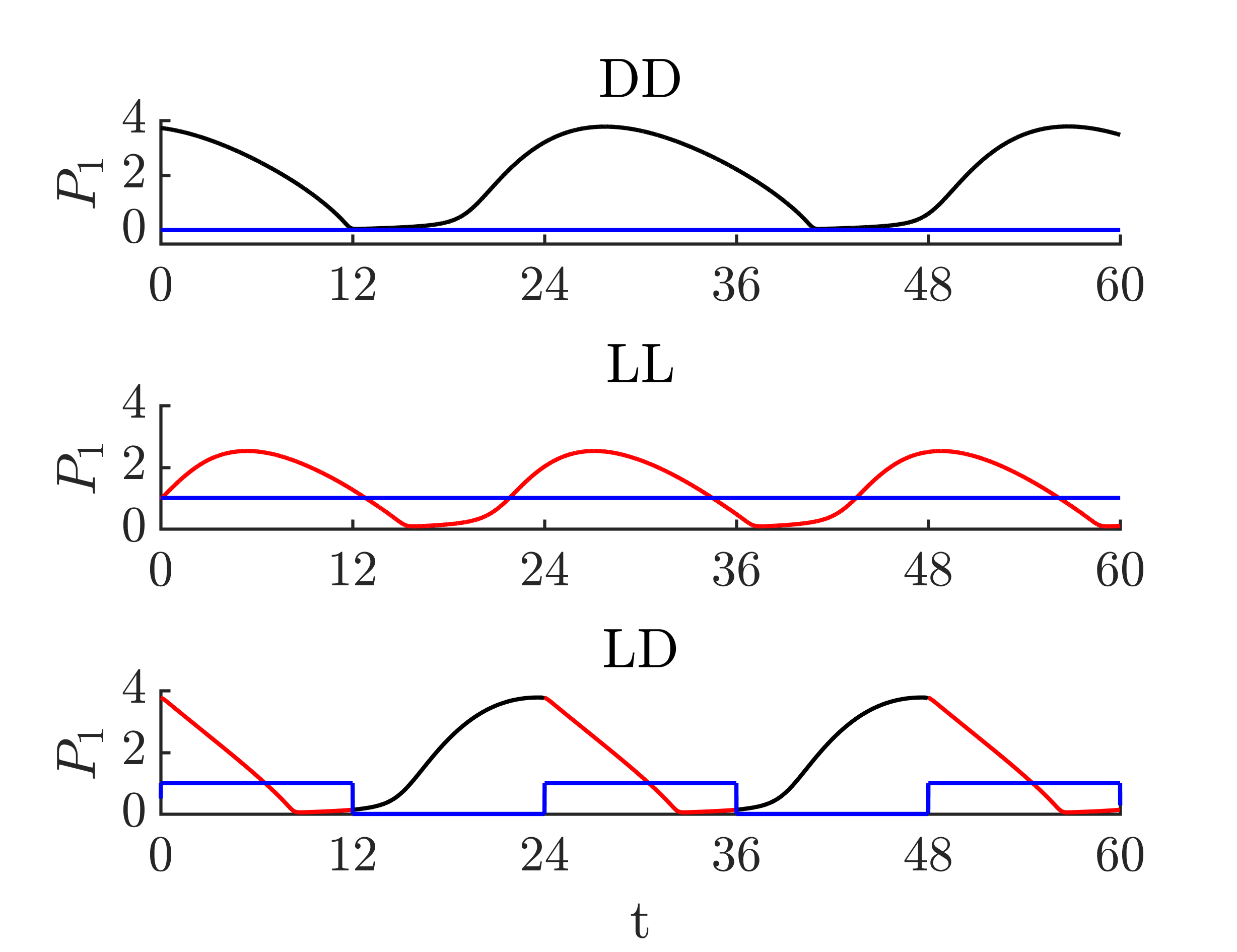}
			\caption{}
			\label{fig:O1_timecourse}
		\end{subfigure}
		\begin{subfigure}{.45\linewidth}
			\centering
			\includegraphics[width=\linewidth]{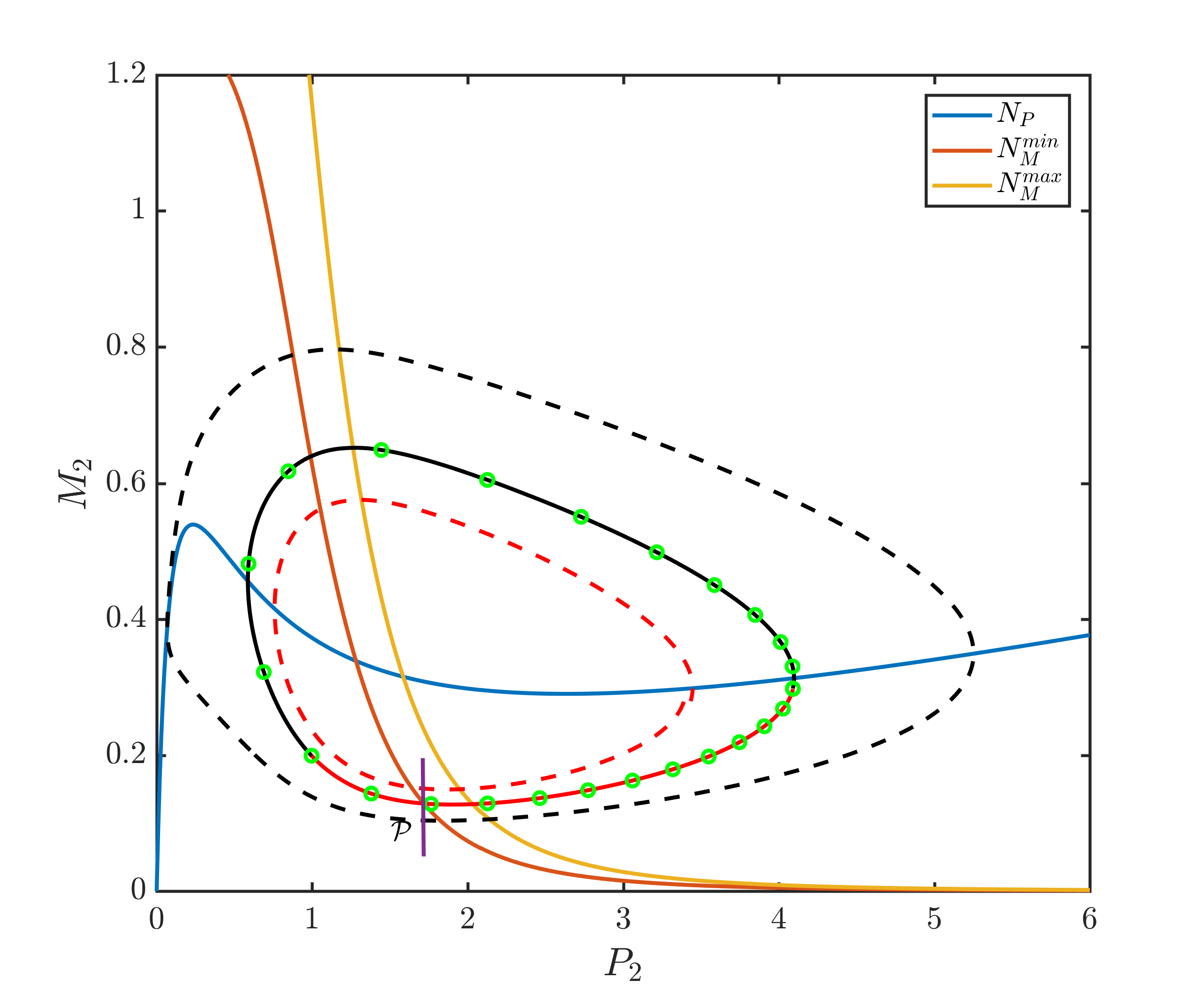}
			\caption{}
			\label{fig:O2_limitcycle}
		\end{subfigure}
		\begin{subfigure}{.45\linewidth}
			\centering
			\includegraphics[width=\linewidth]{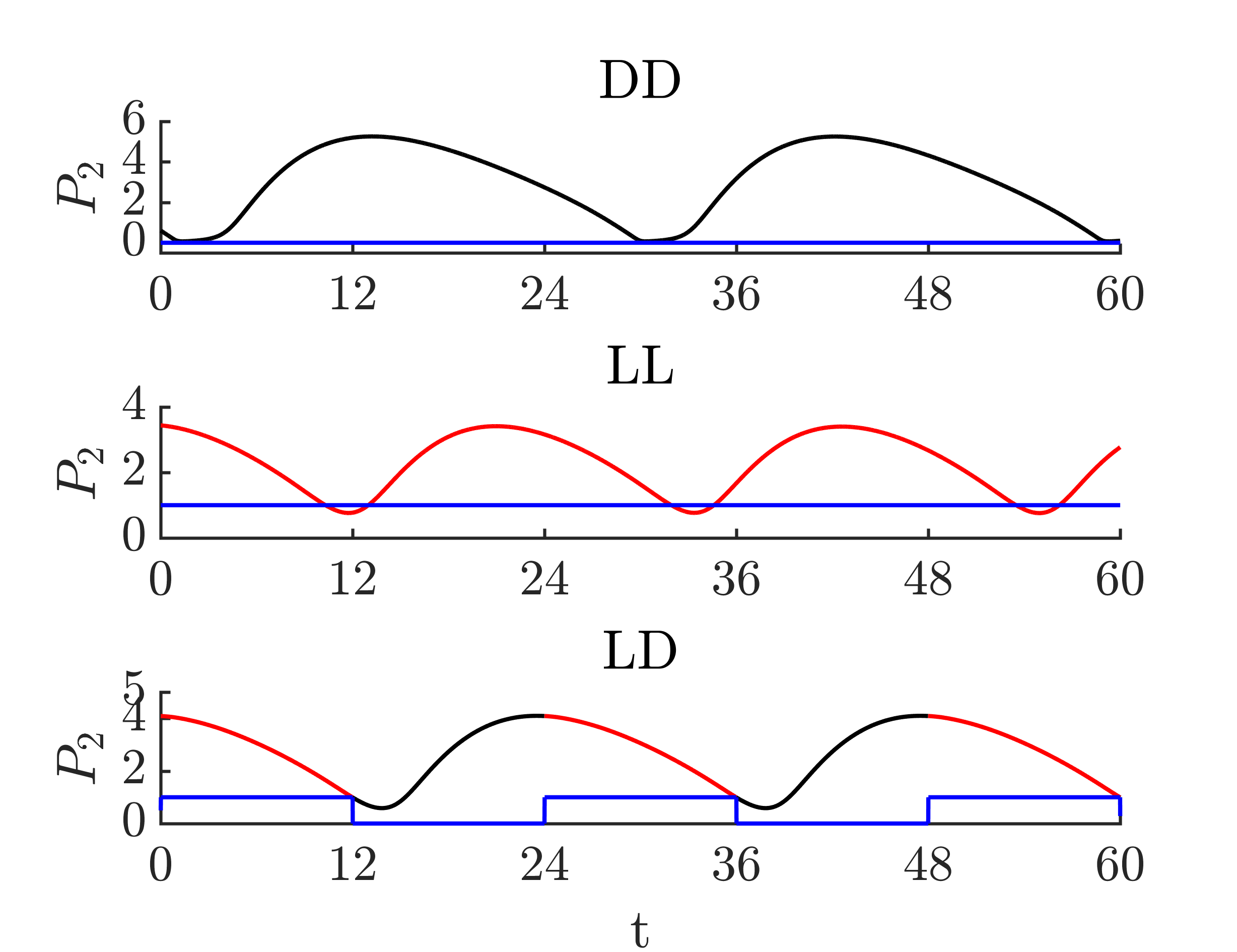}
			\caption{}
			\label{fig:O2_timecourse}
		\end{subfigure}
		
		\caption{(a): The periodic solutions of $O_1$ in DD, LL and LD conditions. The dashed black trajectory represents the DD limit cycle, the dashed red trajectory represents the LL limit cycle. The solid trajectory represents the LD solution with green hourly markers. (b): The time course plots: $P_1$ vs $t$ in all three cases (blue line lies at 0, 1 or is a square wave for DD, LL or LD, respectively.  (c): The periodic solutions of $O_2$ when $O_1$ is in DD, LL and LD conditions. Same color scheme as in (a). The Poincar\'e section is represented at $P_2=1.72$ by a small vertical line segment. (d): The time course plots: $P_2$ vs $t$ in DD, LL and LD conditions.}
		\label{fig:LD}
	\end{figure}

	We plot the entrained solution of the CNT by direct simulation.  In our simulations, we take a specific set of parameters for equation (\ref{CNT}), i.e. $\phi_1=\phi_2=2.1$, $\epsilon_1=\epsilon_2=0.05$, $k_D=0.05,\ k_{L_1}=0.05,\ k_{L_2}=0,\ k_f=1,\ \alpha_1=2$. In Fig.~\ref{fig:O1_limitcycle},  the periodic solutions of $O_1$ are presented for different light conditions. The dashed black (red) limit cycle denotes the stable solution of $O_1$ in DD (LL) conditions.  The solid red-black limit cycle denotes the LD-entrained solution of $O_1$, with hourly markings shown by green open circles.  We also show various nullclines and note that the $M$ nullcline (yellow curve) is unique, but the $P$-nullcline (red and blue curve) varies between $(k_D+k_L)P + k_f h(P)$ and $k_D P + k_f h(P)$. The corresponding time courses are shown for the $P_1$ variable in Fig. \ref{fig:O1_timecourse}.
	
	In Fig.~\ref{fig:O2_limitcycle}, we show the entrained solutions of $O_2$ when $O_1$ is in different light conditions. The color convention is the same as in Fig.~\ref{fig:O1_limitcycle}. Here, we note that the $P$-nullcline (blue curve) is unique, but the $M$-nullcline (red and yellow curves) varies between $(1+\alpha_1 min|M_1(t)|)g(P_2)$ and $(1+\alpha_1 max|M_1(t)|)g(P_2)$. The Poincar\'e section is represented as a small purple vertical line segment at $P_2=1.72$. For completeness, we also show the time course plots related to the same condition in Fig.
	\ref{fig:O2_timecourse}. The time course plots show that the period of the DD solution is longer than that of LD, and the period of LL solution is shorter than that of LD. In particular, we found that the period of the DD cycle is 28.9 h, which is the same as the DD cycle of $O_1$, and the period of the LL cycle is 21.6 h, which is also the same as the LL cycle of $O_1$. This is not surprising, because when the coupling strength is strong enough, $O_2$ is entrained by $O_1$.\\

Remark: Justification of the choice of Poincar\'e section\\
		We claim that $\mathcal{P}$ is a global section. Which means taking any $(P_2,M_2) \in \mathcal{P}$, it will eventually return to $\mathcal{P}$.\\
		Based on phase plane analysis of Fig. \ref{fig:O2_limitcycle}, it is easy to show that any trajectory starting on $\mathcal{P}$ will cross the right branch of $N_P$, we will show that those trajectories stay close in $M$ direction.\\
		Suppose we have a trajectory cross the right branch of $N_P$ at $(\tilde{P}_2,\tilde{M}_2)$, where $\tilde{P}_2>3$, so that
		\begin{equation*}
		\tilde{M}_2 = H(P_2) = k_{f} h(\tilde{P}_2) + k_{D}\tilde{P}_2
		\end{equation*}
		Take a derivative of the right hand side function, for convenient, we use $x$ repsents the $P_2$ variable,
		\begin{equation*}
		H'(x) = k_f h'(x) + k_D = k_f \frac{0.1-x^2}{(0.1+x+x^2)^2} + k_D
		\end{equation*}
		When $x$ is large, $H'(x) \rightarrow k_D$, $H(x) \approx k_D x$, and $k_D$ is selected to be a small parameter, so when $\tilde{P}_2>3$,
		\begin{equation*}
		H(x_1) - H(x_2) \approx k_D(x_1-x_2)
		\end{equation*} 
		So the difference of $M_2$ between two points on the right branch of $N_P$ is small. Next we show that those points have approximately the same dynamics along $M_2$ direction.
		When $P_2$ is large, $g(P_2)\rightarrow 0$, the second equation of \eqref{Reduced_CNT} is approximately
		\begin{align*}
		\frac{dM_2}{dt} &=  -\phi_2 \epsilon M_2 \\
		M_2(t) 				   &= C_0 e^{-\phi_2\epsilon t}
		\end{align*}
		The main point here is that the effect of $M_1$ is gone, so the solution is continuously dependent on initial condition. We already showed that the initial points $(\tilde{P}_2,\tilde{M}_2)$ is close on $M_2$ direction. So we can argue that the they must stay close in $M_2$ direction when they hit the Poincar\'e section.

	\subsection{The $O_1$-entrained map}
	
\begin{figure}[H]
		\centering
		\begin{subfigure}{0.45\linewidth}
			\includegraphics[width=\linewidth]{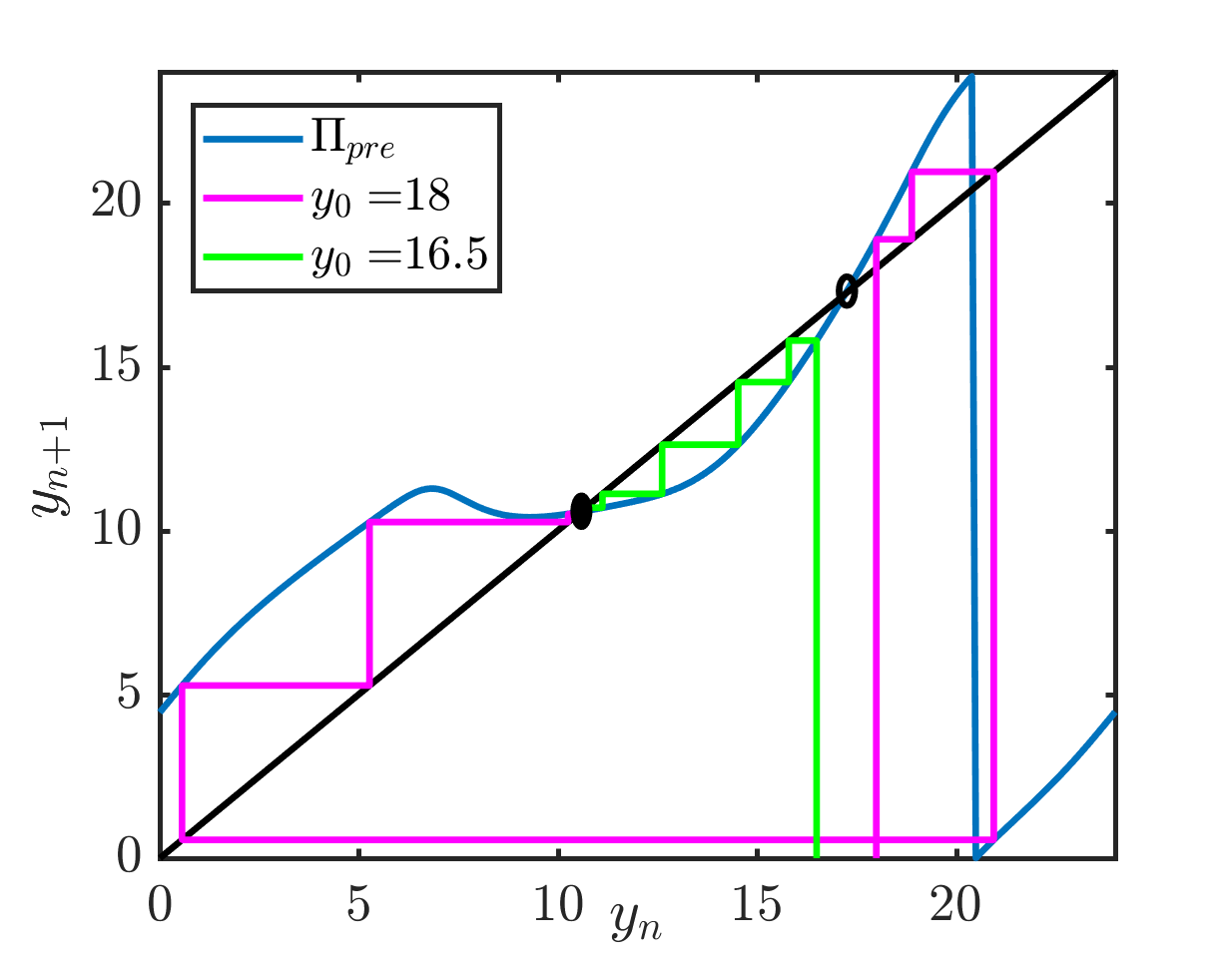}
			\caption{}
		\end{subfigure}
		\begin{subfigure}{0.45\linewidth}
			\includegraphics[width=\linewidth]{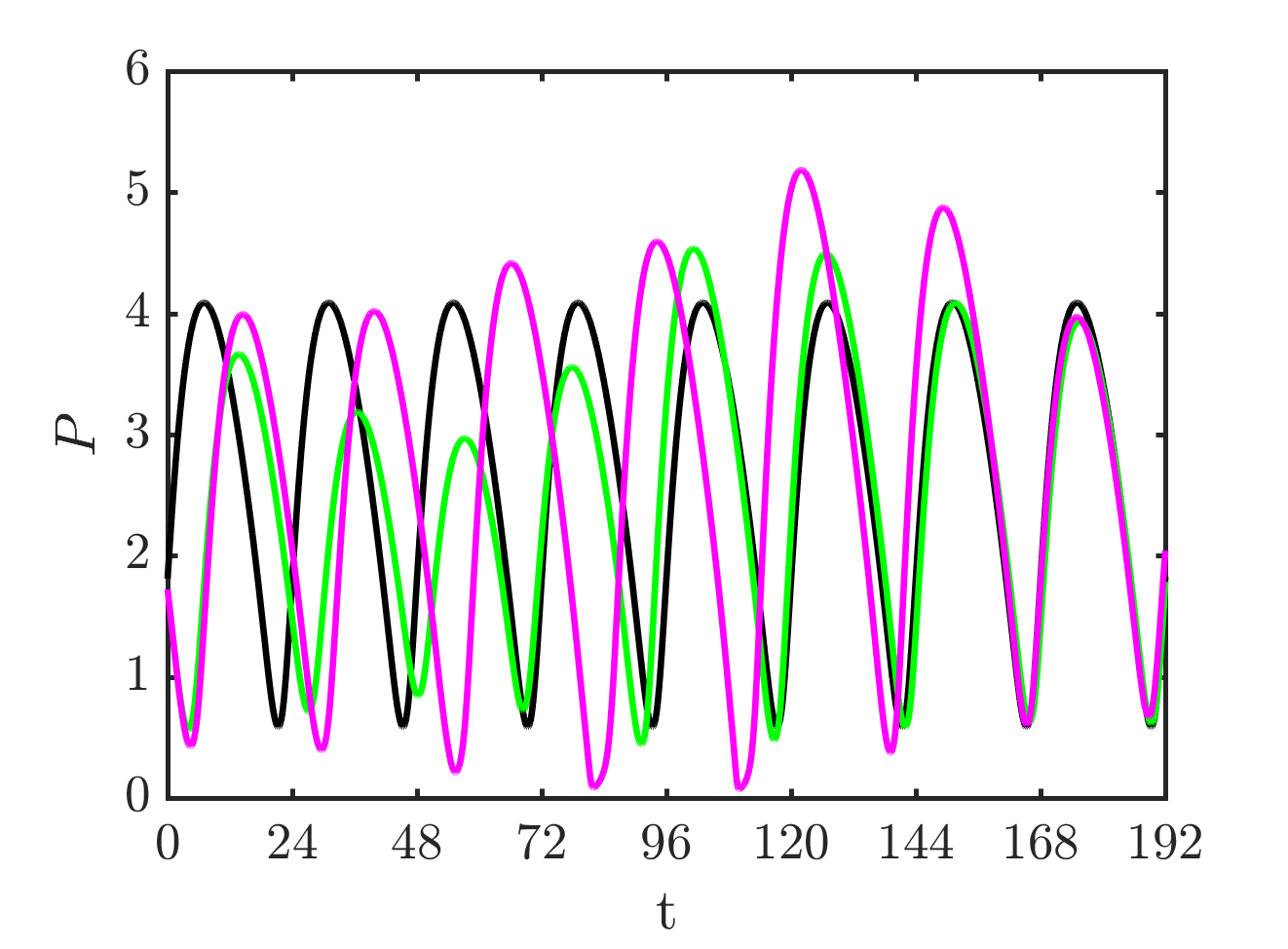}
			\caption{}
		\end{subfigure}
		\begin{subfigure}{0.45\linewidth}
			\includegraphics[width=\linewidth]{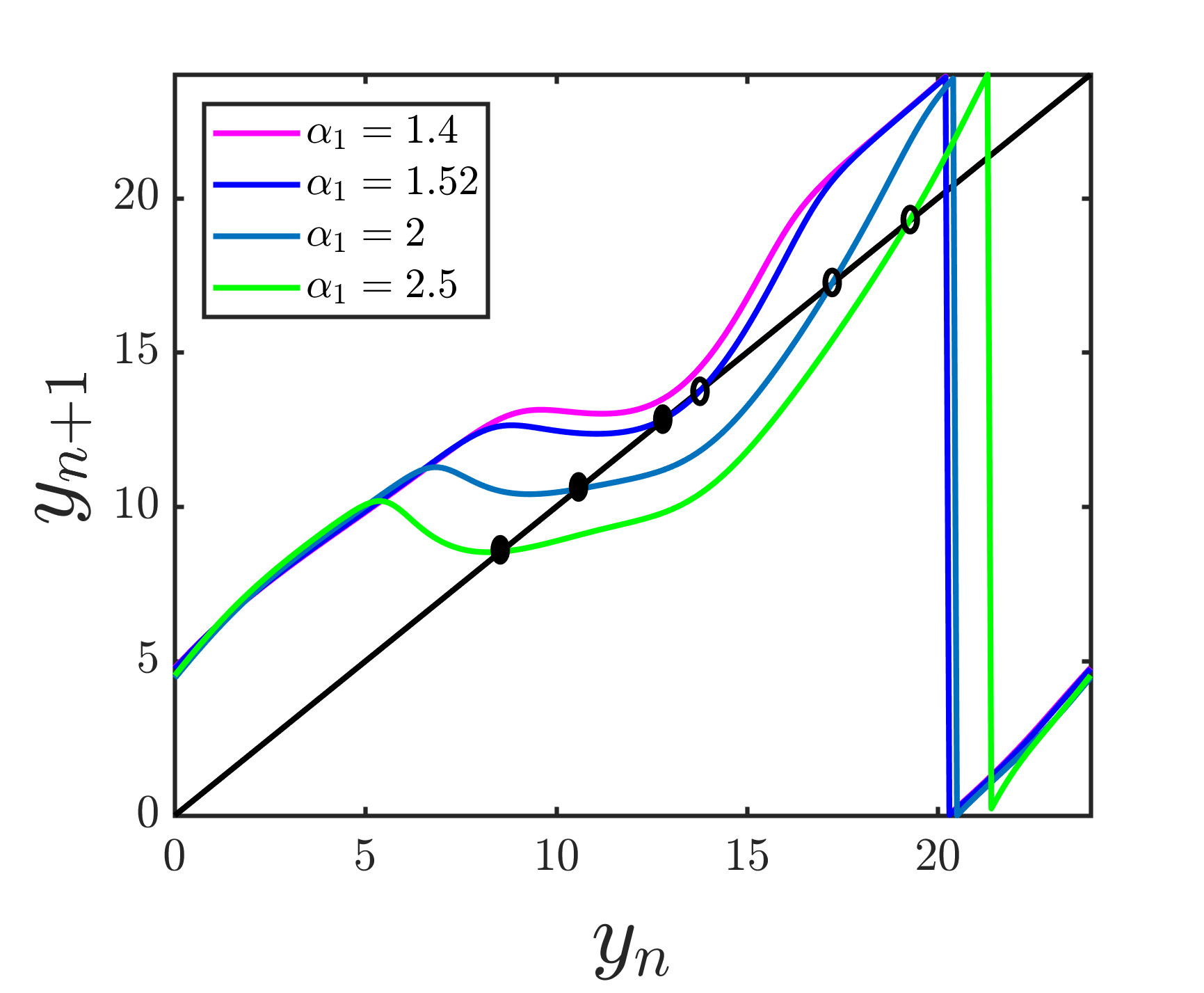}
			\caption{}
		\end{subfigure}
		\begin{subfigure}{0.45\linewidth}
			\includegraphics[width=\linewidth]{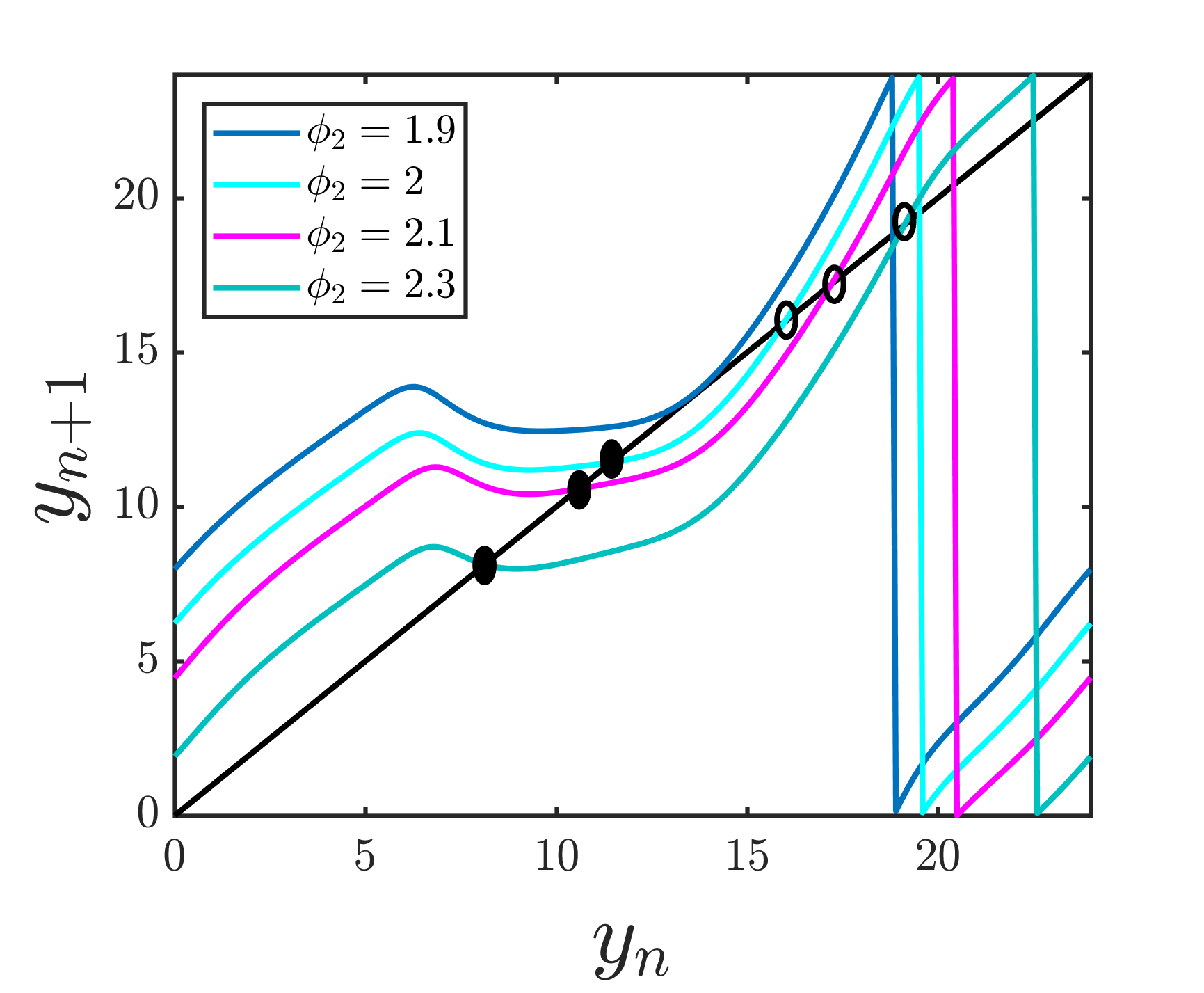}
			\caption{}
		\end{subfigure}
		\caption{(a): The cobweb diagram for the $O_1$-entrained map. We pick two different initial conditions and show how the iterates move to the stable fixed point. (b): The approach to the stable solution (black curve) in the $t$ vs $P$ plane; the colors correspond to the two initial conditions in (a). (c): The map displays a saddle-node bifurcation by decreasing $\alpha_1$. (d): Decreasing the intrinsic period of $O_2$ by decreasing $\phi_2$ also leads the map to display a saddle-node bifurcation. Fixed points shown as open circles are unstable, and those shown with solid circles are stable.}
		\label{fig:1DPEMap}
	\end{figure}

	The $O_1$-entrained map we obtained from Eq. \eqref{Eq:premap} has similar properties as the entrainment map Diekman and Bose constructed in their paper \cite{diekman2016entrainment}. Figure \ref{fig:1DPEMap} shows that there are two fixed points which correspond to different types of periodic solutions for the CNT system. The lower one with $y_{n+1}=y_n=10.2$ is a stable fixed point of the map, which represents a stable periodic solution. The upper one with $y_{n+1}=y_n=17.2$ is an unstable fixed point of the map.
	
	We classify the direction of entrainment as occurring through phase advance or phase delay. Suppose $y_{n+1}=\Pi_{O_1}(y_n)$, and the return time needed from $y_n$ to $y_{n+1}$ is less than 24 hours. We call this a phase advance. Alternatively, if the return time is greater than 24 hours, we call it phase delay.  The unstable fixed point of the map plays an important role in determining this direction. For example, pick two different initial conditions ($y_0=16.5,18$) near the unstable fixed point and  use the cobweb method to observe how different directions of entrainment can occur. For $y_0=16.5$, the iterates move to the left and converge to the stable solution by phase advance. For $y_0=18$ however, the iterates move to the right and converge to the stable solution by phase delay. In Fig. \ref{fig:1DPEMap}(b), we compare the iterates with simulations; the green curve corresponds to $y_0=16.5$ and the magenta curve corresponds to $y_0=18$. The black curve is the entrained solution for $O_2$. The direction of entrainment from the simulations agrees with the  calculations obtained from the map.
	
	In our model system, there are two parameters of interest, the coupling strength $\alpha_1$ and the intrinsic period of $O_2$ governed by $\phi_2$. In Fig. \ref{fig:1DPEMap}(c), we decrease $\alpha_1$ from 2.5 to 1.4, so that the coupling strength is weaker. As a result, the return time $\rho(y)$ increases. This makes the map move up, and the stable and unstable fixed points get closer to each other. At $\alpha_1=1.51$, the two fixed points collide at a saddle-node bifurcation. In Fig. \ref{fig:1DPEMap}(d), we increase the intrinsic period of $O_2$ by decreasing $\phi_2$ from 2.3 to 1.9, so that the difference between the intrinsic period and the 24-h forcing increases, which increases the return time to the Poincar\'e section. Hence the map moves up. When $\phi_2=1.91$, the map passes through the saddle-node bifurcation value. Notice that, the fixed point of the map corresponds to a 1:1 phase locked solution of the full system. When we lose the fixed point in the map, which means we also lose the entrained solution of the full system, and we did some simulations with the bifurcation value of parameters $\alpha_1 \approx 1.51$, $\phi_2 \approx 1.91$ of the full system. We do see the lose of entrainment.
	
\begin{figure}[H]
		\centering
		\begin{subfigure}{0.45\linewidth}
			\centering
			\includegraphics[width=\linewidth]{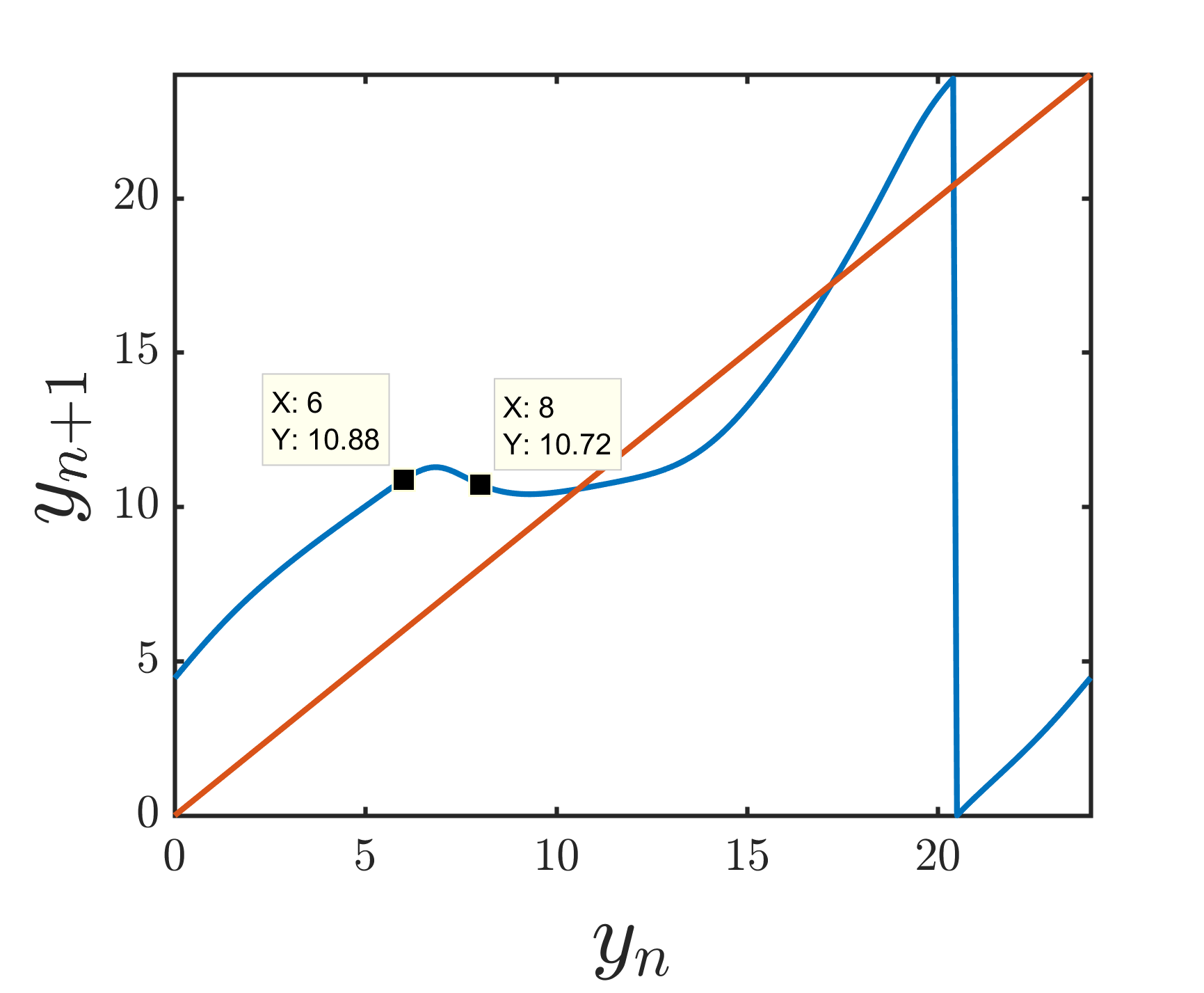}
			\caption{}
		\end{subfigure}
		\begin{subfigure}{0.45\linewidth}
			\centering
			\includegraphics[width=\linewidth]{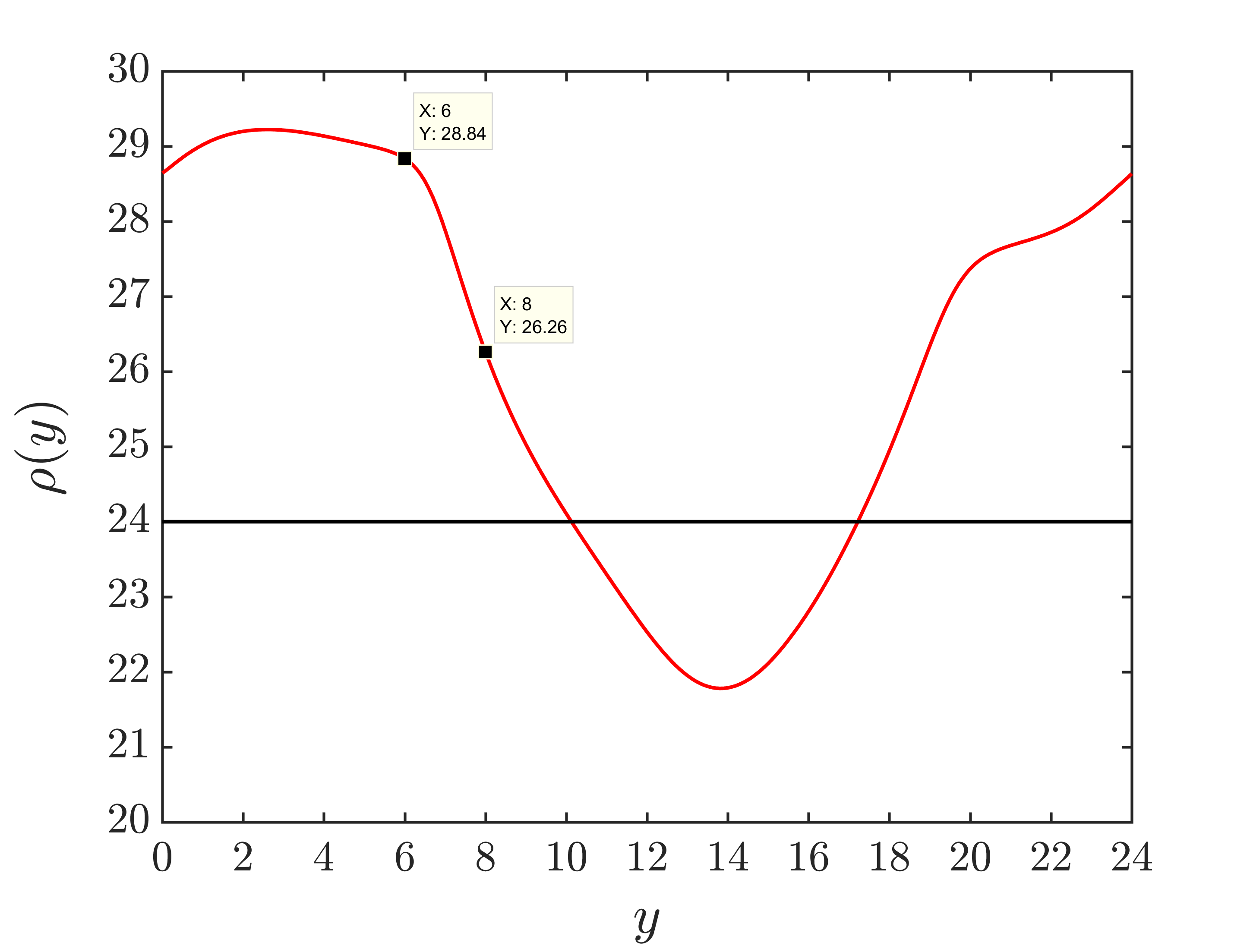}
			\caption{}
		\end{subfigure}
		\begin{subfigure}{.45\linewidth}
			\centering
			\includegraphics[width=\linewidth]{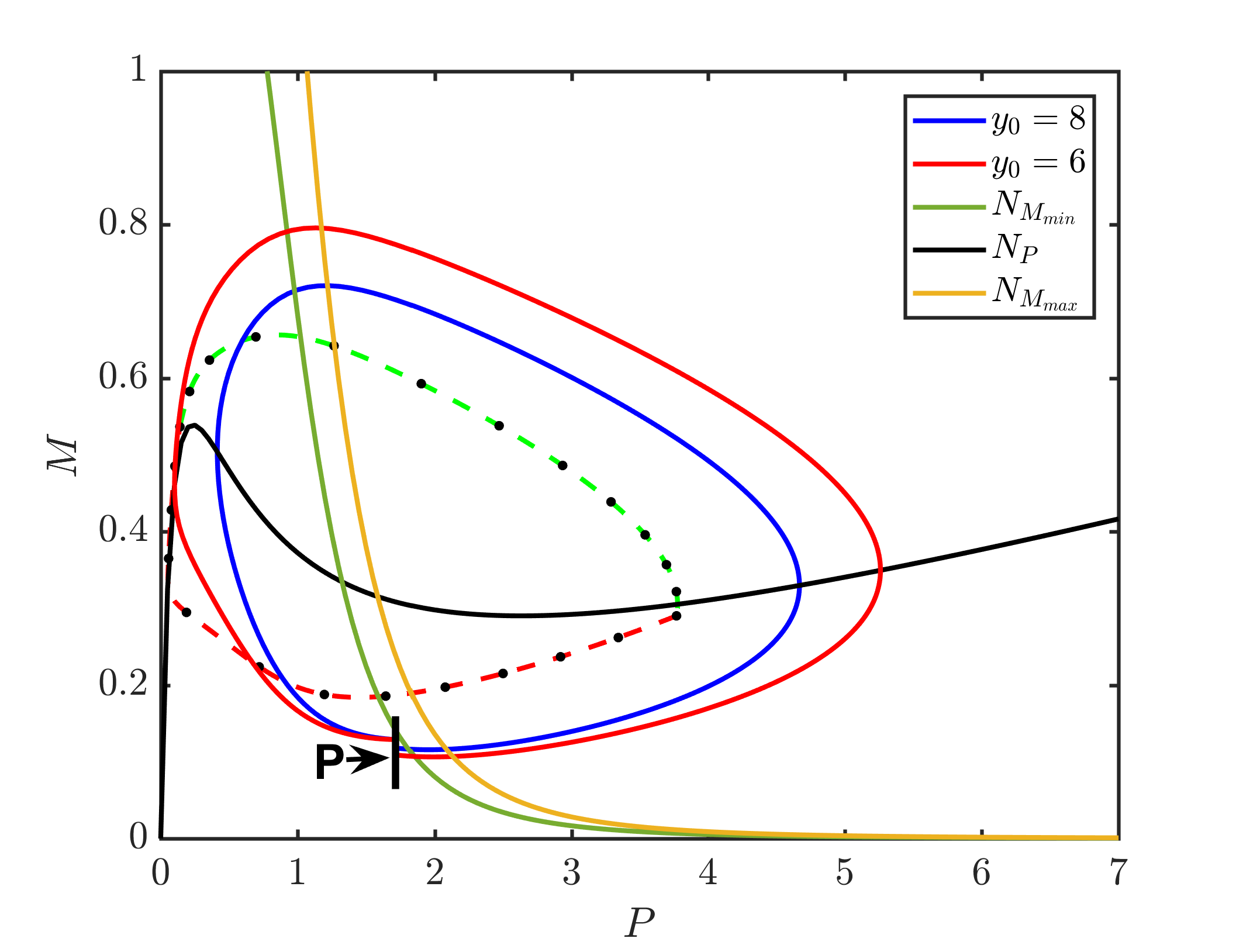}
			\caption{}
		\end{subfigure}
		\begin{subfigure}{.45\linewidth}
			\centering
			\includegraphics[width=\linewidth]{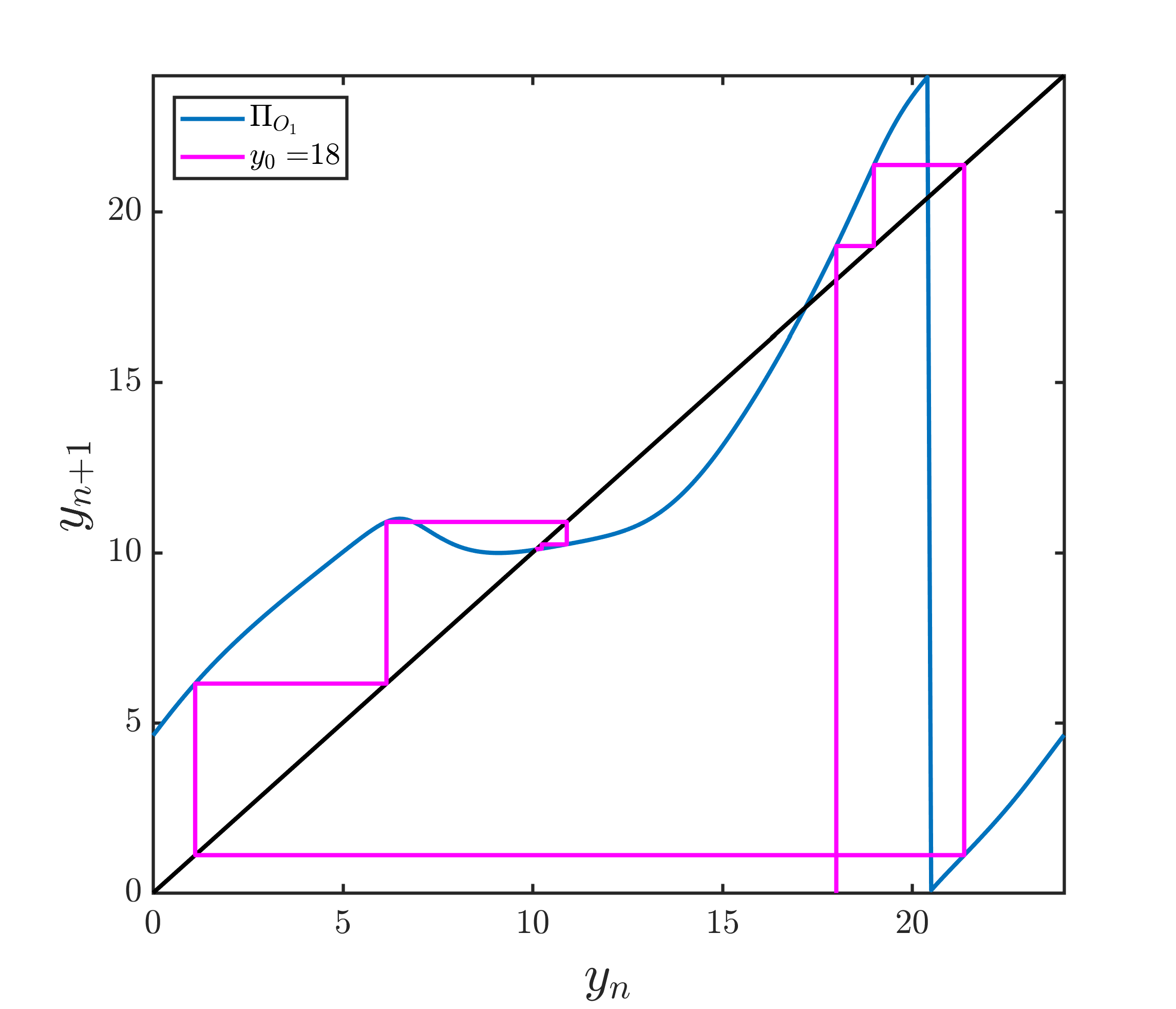}
			\caption{}
		\end{subfigure}
		
		\caption{Non-monotonicity in the entrainment map leads to convergence initially due to phase delay but ultimately due to phase advance.  (a): The non-monotone $O_1$-entrained map and two choices of initial conditions near the local maximum. Note that the local maximum lies above the value of the fixed point of the map. (b): The return time plot associated with the two initial conditions. (c): The corresponding phase plane. The solid blue trajectory for $y_0=8$ does not approach the left branch of $N_P$, while the solid red trajectory for $y_0=6$ does, causing its evolution to slow down. (d): Starting with an initial condition $y_0=18$, the first four iterates phase delay. The fourth iterate lands near the local max of the map, and subsequent iterates then phase advance.}
		\label{fig:nonmono}
	\end{figure}

	Notice that the $O_1$-entrained map we construct is not monotonic, which makes it different from the 1-D entrainment map found in \cite{diekman2016entrainment}. To understand this nonmonotonicity, we take two initial conditions ($y_0=6$ and $y_0=8$) near the local maximum of the map in Fig. \ref{fig:nonmono}(a), and analyze the dynamics of the system. Associated with the return time plot in Fig. \ref{fig:nonmono}(b), we found that the return time is between 28 and 29 when $y$ is less than the local maximum point. But when it crosses that point, the return time decreases quickly with the derivative $\rho'(y)<-1$. In Fig. \ref{fig:nonmono}(c), we plot the trajectories with the two initial conditions. The trajectory for $y_0=6$ flows to the left branch of the $P$-nullcline, which increases the return time since evolution near this branch is slow. Alternatively, the
	trajectory for $y_0=8$ doesn't flow near the left branch and thus has a shorter return time. A minor consequence of this non-monotonicity is that some solutions converge to the stable fixed point by initially phase delaying, but then ultimately phase advancing. For example, in Fig. \ref{fig:nonmono}(d), we take $y_0=18$ then cobweb the map. We find that the first four iterates initially phase delay. The fourth iterate lands near the local maximum of the map, which lies above the value of the fixed point. This causes subsequent iterates to phase advance. This non-monotonicity foreshadows a more complicated picture that arises under the dynamics of the 2-D map.

	\subsection{The results of the general 2-D map}
	
	In this section, the analysis of the 2-D map is presented. We follow ideas first derived by Akcay et. al. \cite{akcay2014effects} and followed up on in \cite{akcay2018phase} to find  fixed points of the map via a geometric method. The entrainment time and the direction of entrainment are analyzed by iterating the map. We also compare these results with simulations. At the end of this section, we show that the map is also applicable to the semi-hierarchical model.
	
\begin{figure}[H]
		
		\centering
		\begin{subfigure}{.45\linewidth}
			\centering
			\includegraphics[width=\linewidth]{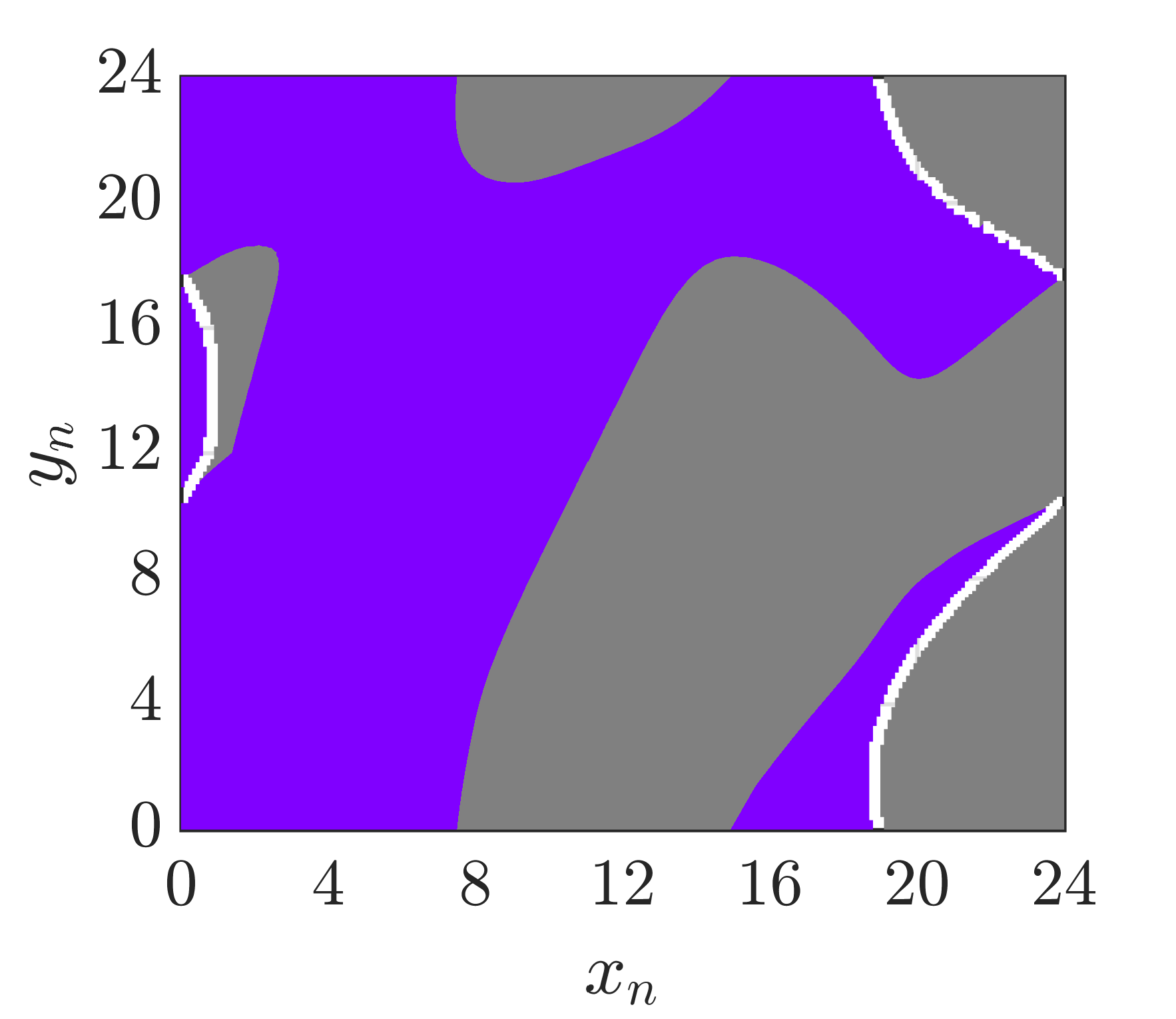}
			\caption{}
			\label{fig:surf_pi1_top}
		\end{subfigure}
		\begin{subfigure}{.45\linewidth}
			\centering
			\includegraphics[width=\linewidth]{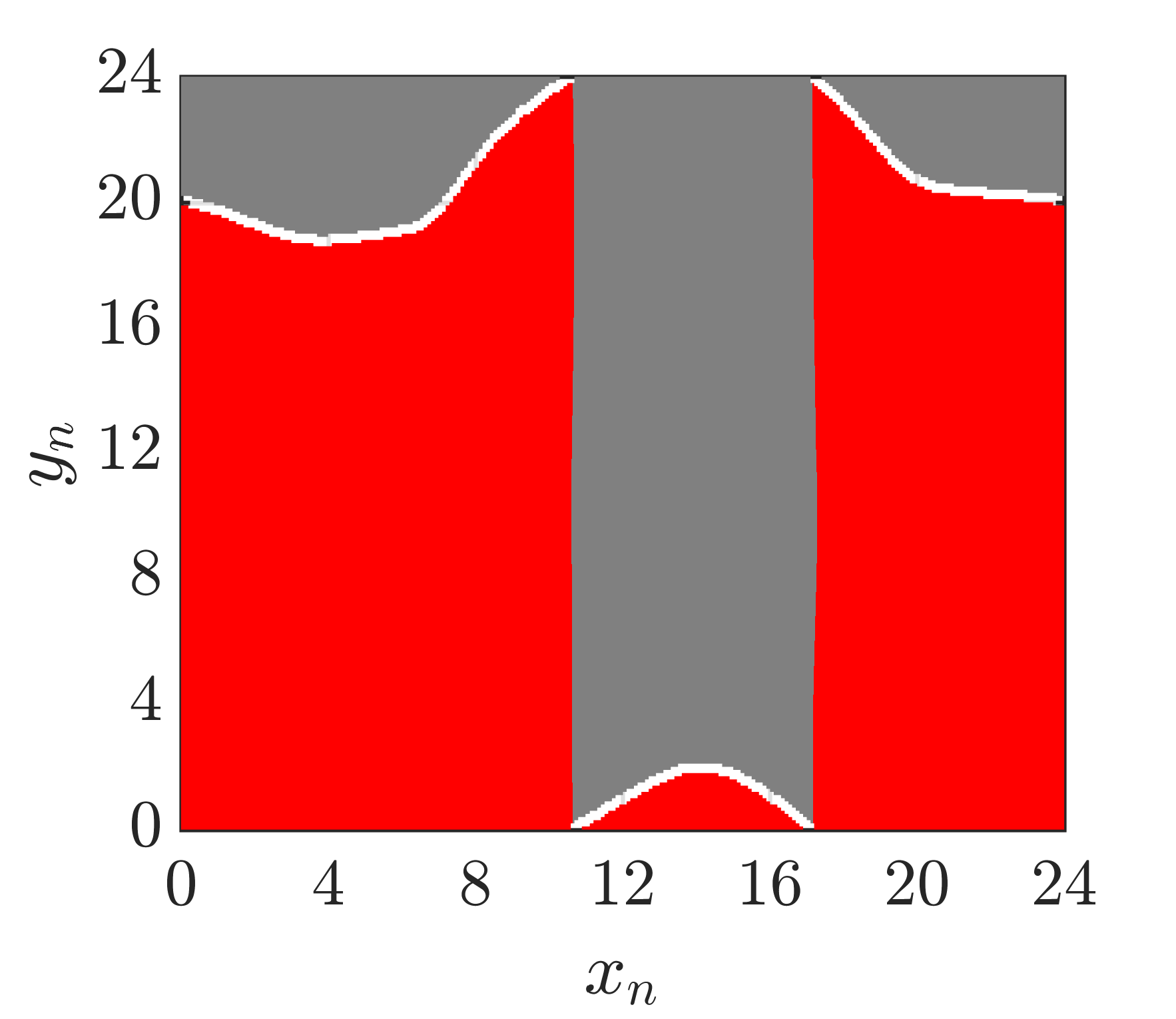}
			\caption{}
			\label{fig:surf_pi2_top}
		\end{subfigure}
		\begin{subfigure}{.45\linewidth}
			\centering
			\includegraphics[width=\linewidth]{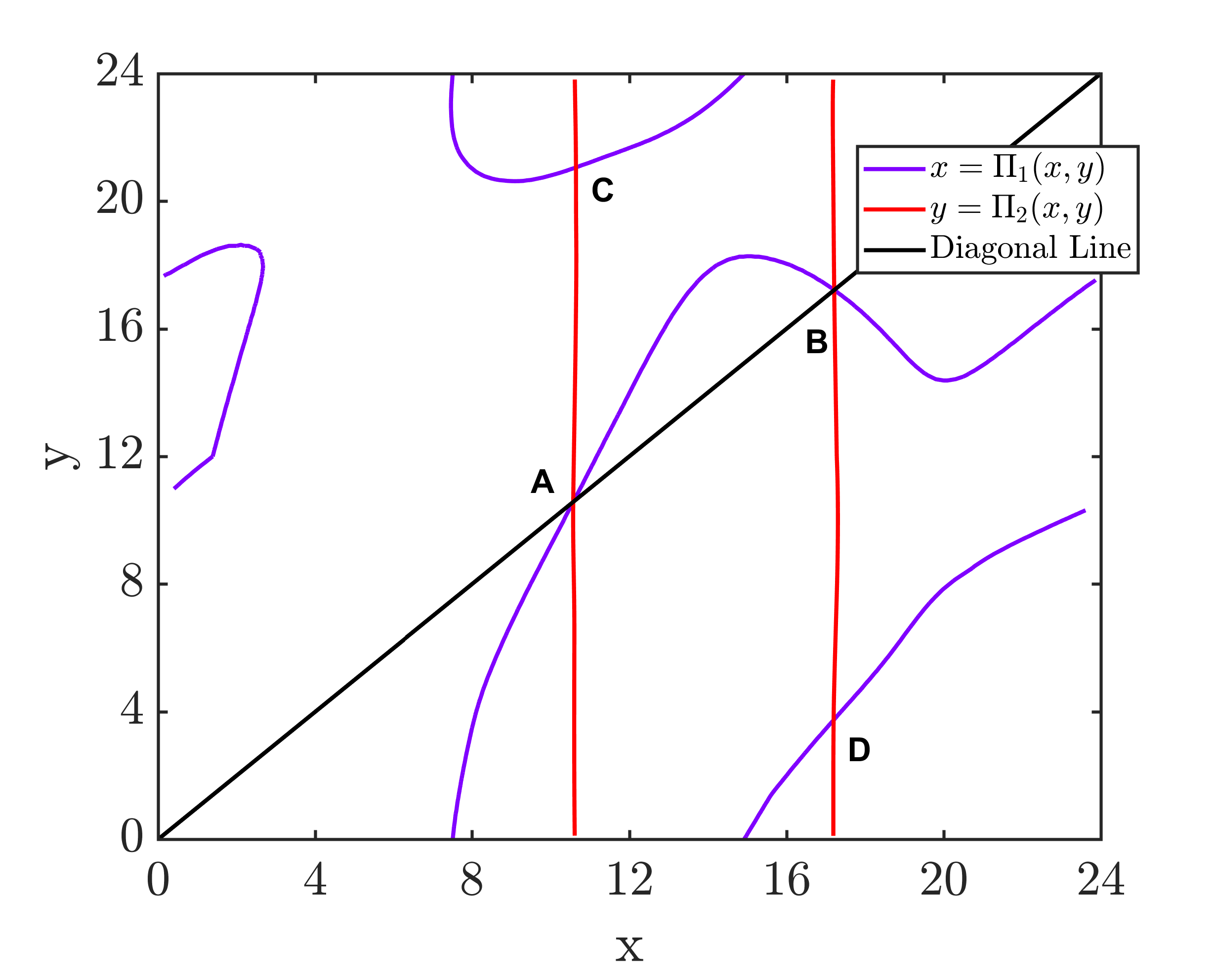}
			\caption{}
			\label{fig:contours}
		\end{subfigure}
		\begin{subfigure}{.45\linewidth}
			\centering
			\includegraphics[width=\linewidth]{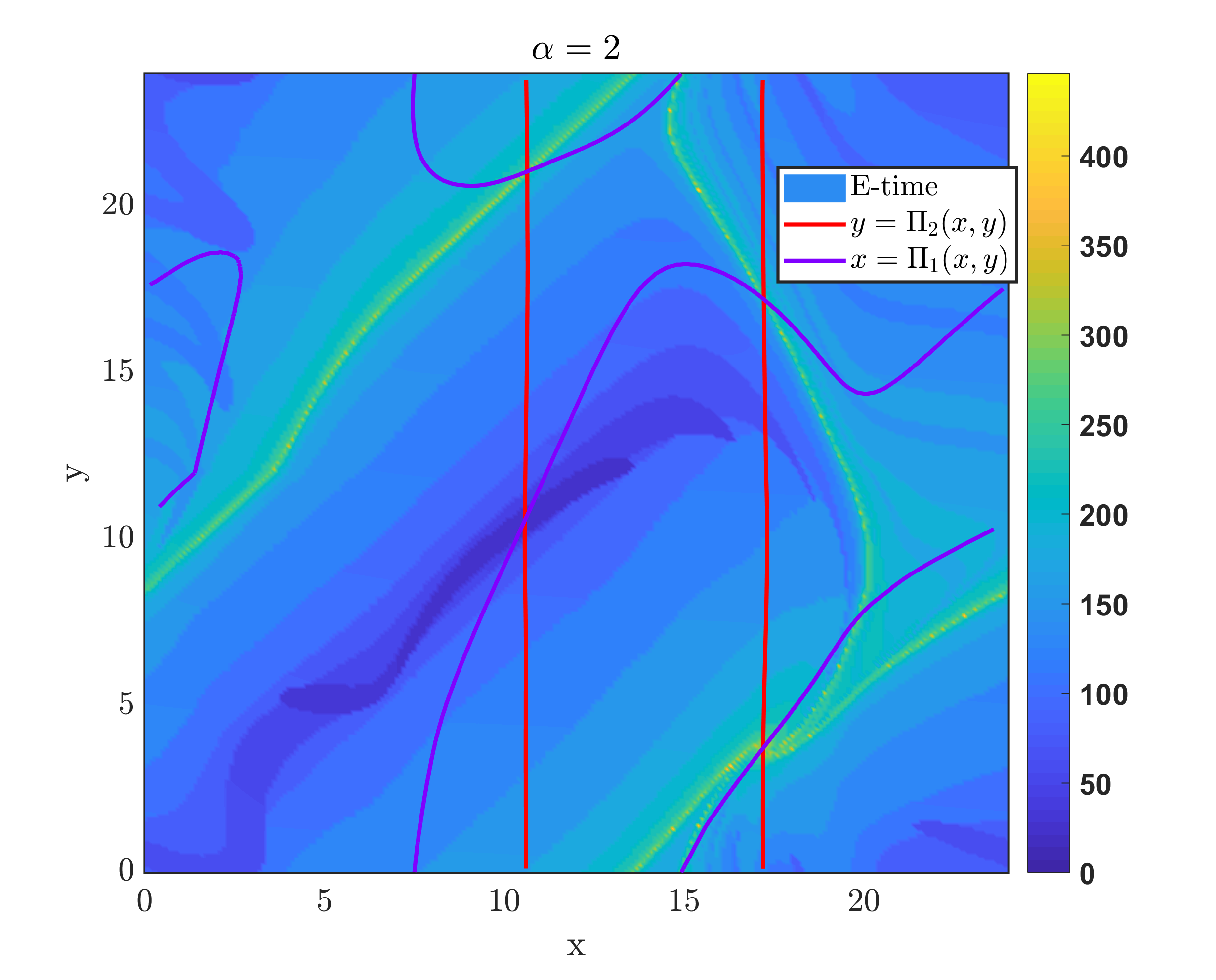}
			\caption{}
			\label{fig:Etimemap_alpha_2}
		\end{subfigure}
		\caption{(a) and (b): The 2-D entrainment map is plotted as two separate maps $\Pi_1\text{ and }\Pi_2$, and projected onto the domain space $(x_n,y_n)$. The purple and red color in both maps denote all points that are above the diagonal plane. The grey color denotes points that are below the diagonal plane. The white curves denote the discontinuity. (c): The purple curves denote points of $\Pi_1$'s nullcline $N_x$ where $x=\Pi_1(x,y)$, the red curves denote points of $\Pi_2$'s nullcline  $N_y$ where $y=\Pi_2(x,y)$. 	Their intersections are the four fixed points of the map. (d): The entrainment time is plotted with a heatmap. The color denotes the entrainment time starting from a specific initial condition. The light green curves locate $W^s(B)$ and $W^s(C)$ from near which the longest entrainment times occur.}
	\end{figure}

	\paragraph{Basic results from the map}
	Both parts of the 2-D map $\Pi_1\mbox{ and }\Pi_2$ are surfaces in relevant 3-D spaces. Because of the mod 24 operation, each surface will contain discontinuities. In Fig. \ref{fig:surf_pi1_top} and \ref{fig:surf_pi2_top}, we project the surface onto the $x-y$ plane. For $\Pi_1$, the purple part of the surface are points lying above the diagonal plane $z=x$,  in other words, $x_{n+1}>x_n$. The red part of the surface of $\Pi_2$ are points lying above the diagonal plane $z=y$, i.e. $y_{n+1}>y_n$.
	The points of grey color denote all points that are below the diagonal planes, $x_{n+1} < x_n$ and $y_{n+1} <y_n$. The white curves indicate locations of discontinuity of the map. The separation of the two different colors are curves which indicate the points where $x=\Pi_1(x,y)\mbox{ and }y=\Pi_2(x,y)$. Here we define those curves as nullclines of the map:
	\begin{equation*}
	N_x=\{(x,y):x=\Pi_1(x,y)\},\quad N_y=\{(x,y):y=\Pi_2(x,y)\}
	\end{equation*}
	which are plotted in Fig. \ref{fig:contours}. The purple curves denote $N_x$.
	Similarly, the red curves denote $N_y$. Their intersections are four fixed points of the map. We numerically calculated the Jacobian at those fixed points and found the eigenvalues of the linearization. These values and the corresponding stability of each fixed point is shown in Table 1.
	\begin{table}[H]
		\centering
		\begin{tabular}{|l|c|c|c|c|}
			\hline
			& x    & y    & eigenvalue                & stability         \\ \hline
			A & 10.6 & 10.6 & 0.1609,\ \ 0.4453           & sink  \\ \hline
			B & 17.2 & 17.2 & 2.0858,\ \ 0.4238           & saddle \\ \hline
			C & 10.6 & 21.1 & 2.325,\ \ 0.2734           & saddle \\ \hline
			D & 17.2 & 3.7  & 1.595+0.77$i$,\ \ 1.595-0.77$i$ & source \\ \hline
		\end{tabular}
		\caption{Numerical computation of the eigenvalues of the map at the four fixed points. Eigenvalues with modulus less than one correspond to stable directions, while those with modulus greater than one correspond to unstable directions.}
	\end{table}
	From the results of the $O_1$-entrained map, points A and B lying on the diagonal line correspond to the stable solution of $O_1$. For $O_2$, point A corresponds to the stable solution. For point B, the trajectory of $O_2$  returns to the Poincar\'e section after 24 hours but corresponds to the unstable solution of the $O_1$-entrained map. At the fixed point C, $O_1$ lies on its own unstable periodic orbit. This can be inferred from and agrees with the calculation of Diekman and Bose \cite{diekman2016entrainment} who showed that the original 1-D entrainment map has an unstable fixed point that corresponds to an unstable periodic orbit.  Thus $O_1$ is entrained to a 24-hour LD cycle and provides a 24-hour forcing to $O_2$. From simulation, we found that the trajectory of $O_2$ stays for several cycles near what appears to be a stable limit cycle, though it is different from the limit cycle corresponding to point A since $O_1$ is unstable and the forcing signal to $O_2$ is different. At point D, if we check the difference between C and D, we can see that
	\begin{equation*}
	(x_D,y_D)=(x_C,y_C)+6.6\mbox{ mod }24
	\end{equation*}
	so $O_1$ is still on its unstable periodic orbit. That is, points C and D represent conditions where the forcing $M_1(t)$ is identical, but just phase shifted by 6.6 hours. Thus $O_2$ still receives 24-hour forcing so we also expect there to exist an unstable $O_2$ limit cycle for this case.
	
	One advantage of the map is its ability to estimate the entrainment time. Starting from different initial conditions, we iterate the map $(x_{n+1},y_{n+1})=\Pi(x_n,y_n)$ until $\|(x_{n+1},y_{n+1})-(x_s,y_s)\|<0.5$, where point A has coordinates $(x_s,y_s)$. Then the entrainment time is the sum of the return times corresponding to each iterate. In Fig. \ref{fig:Etimemap_alpha_2}, we show the entrainment times corresponding to different initial conditions on the torus expanded as a square. We also plot the nullclines $N_x$ and $N_y$  on top of it for illustrative purposes. The color for each point on the square denotes the entrainment time needed for that initial point.

	
	Notice that, in Fig. \ref{fig:Etimemap_alpha_2}, there are two light green curves. Along these curves, the entrainment time is much longer than other regions. Additionally, they appear to connect the two saddle points B, C, with the unstable source D. Though not proven here, we believe that these curves locate where the stable manifolds of the saddle points B and C ($W^s(B)$ and $W^s(C)$) are. To completely understand the dynamics of the entrainment map, it is useful to numerically find the stable and unstable manifolds.
	
	The algorithm we used to find the manifolds of the entrainment map are based on the following results. For the unstable manifold, Krauskopf and Osinga \cite{krauskopf1998growing} introduced a growing method to calculate the unstable manifold point by point. They initially iterate points chosen in a neighborhood of the fixed point along the associated unstable eigenvector and accept new points as lying on the unstable manifold if they satisfy specific constraints. For the stable manifold, the search circle (SC) method introduced by England et. al. \cite{england2004computing} utilizes the stable eigenvector to find points within a certain radius that iterate onto a segment of the stable eigenvector. The SC method has the advantage that it does not require the inverse of the map to exist, which is important for us since our map is non-invertible. Both of these methods are constructed for planar non-periodic domains. In our case, the map lives on a torus, but is graphically shown on a square. Whenever an iterated point exceeds the boundary of the square, we use the modulus operation to define the correct value within the square. Thus we develop our algorithm to account for this discontinuity. Another difference is that the terminating conditions for both the growing and SC methods rely on calculating the arclength of the manifolds up to a certain predetermined length.  However, in our map, the stable manifolds of points B and C are generated from the source point D, while their unstable manifolds terminate at point A. Thus our algorithm terminates when these manifolds enter prescribed  neighborhoods of those corresponding fixed points D and A.
	
\begin{figure}[H]
		
		\centering
		\begin{subfigure}{.45\linewidth}
			\includegraphics[width=\linewidth]{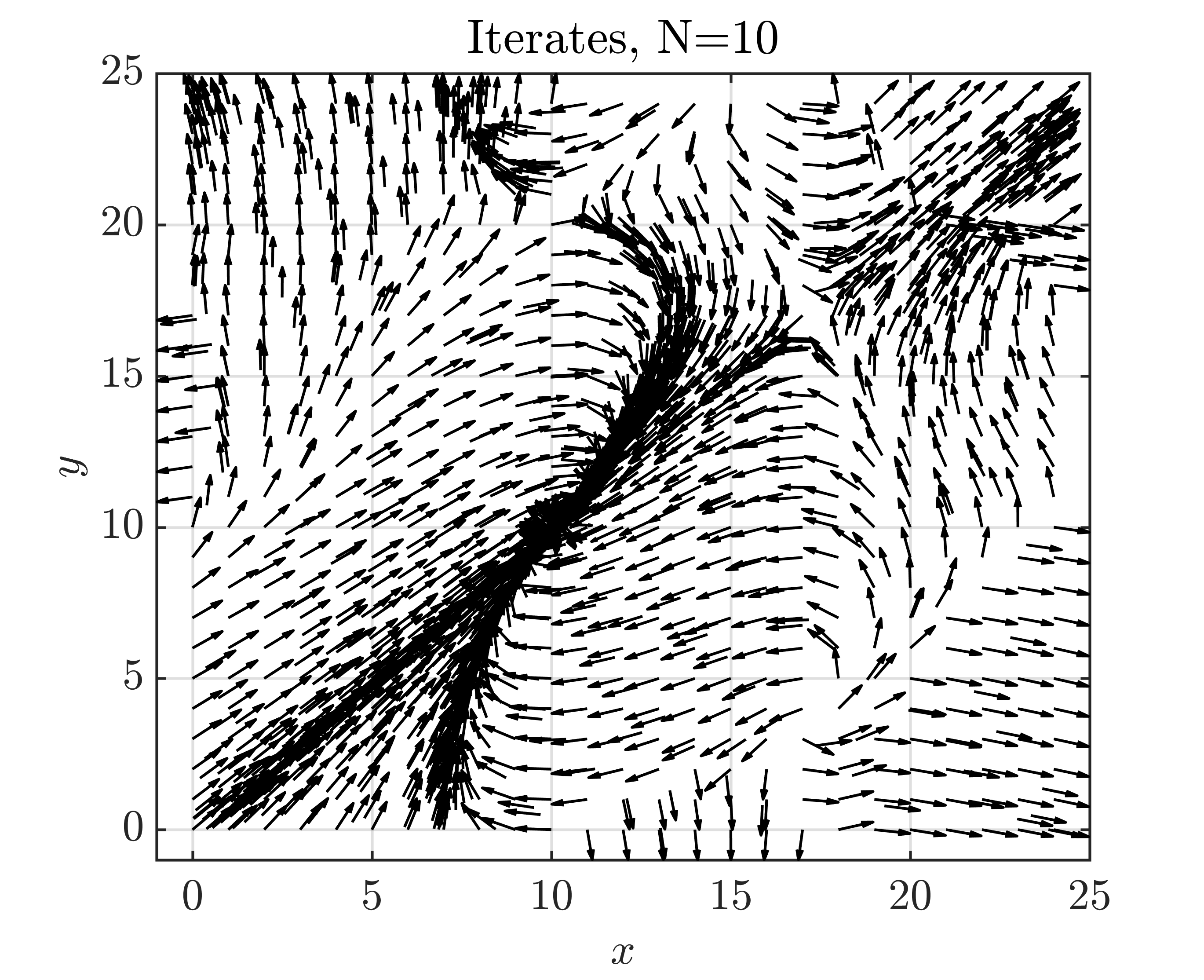}
			\caption{}
		\end{subfigure}
		\begin{subfigure}{.45\linewidth}
			\includegraphics[width=\linewidth]{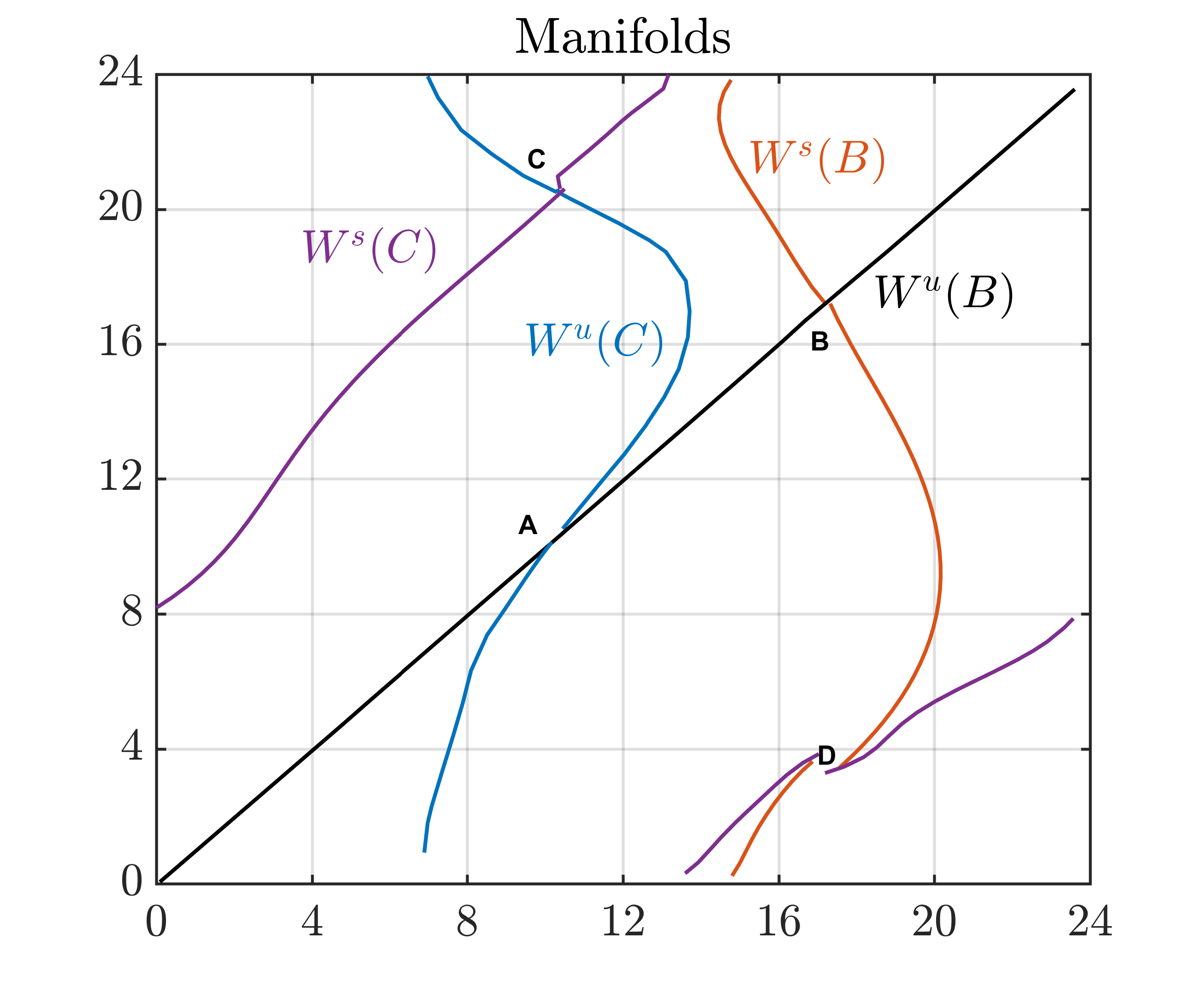}
			\caption{}
		\end{subfigure}
		
		\caption{(a): $N=10$ iterates from various initial points are shown. The arrows at each coordinate point in the direction of the next iterate. The vector field indicates that there may exist a separatrix type structure at both points B and C. (b): Stable and unstable manifolds of B and C as generated through the generalization of the search circle and growing methods (see text). The labeled manifolds do appear to provide a separatrix type behavior despite this being a map and not a flow.}
		\label{fig:vf}
	\end{figure}

	In Fig. \ref{fig:vf}a, we choose initial points ranging from $0<x<24,0<y<24$, and iterate N times for each initial point. The arrows on each coordinate are pointing to its own next iterate. The obtained vector field give us another visualization of the map. In Fig. \ref{fig:vf}b, the numerical result of stable and unstable manifolds of B and C are plotted. $W^s(B)$ and $W^s(C)$  agree with the light green curves in Fig. \ref{fig:Etimemap_alpha_2}. $W^u(B)$ is exactly the diagonal line of the phase plane, which is not surprising. Because the diagonal line corresponds to the $O_1$-entrained case, if an iterate starts on the diagonal line, it stays on it. The numerical calculation of the eigenvector of $E^u(B)$ is approximately (0.7,0.7) on the diagonal line, which means $W^u(B)=E^u(B)$. $W^u(C)$ also matches the darkest region in Fig. \ref{fig:vf}(a). Indeed, these dark regions indicate the location of the unstable manifolds of points B and C.
	
\begin{figure}[H]
		\centering
		\begin{subfigure}{.8\linewidth}
			\includegraphics[width=.5\linewidth]{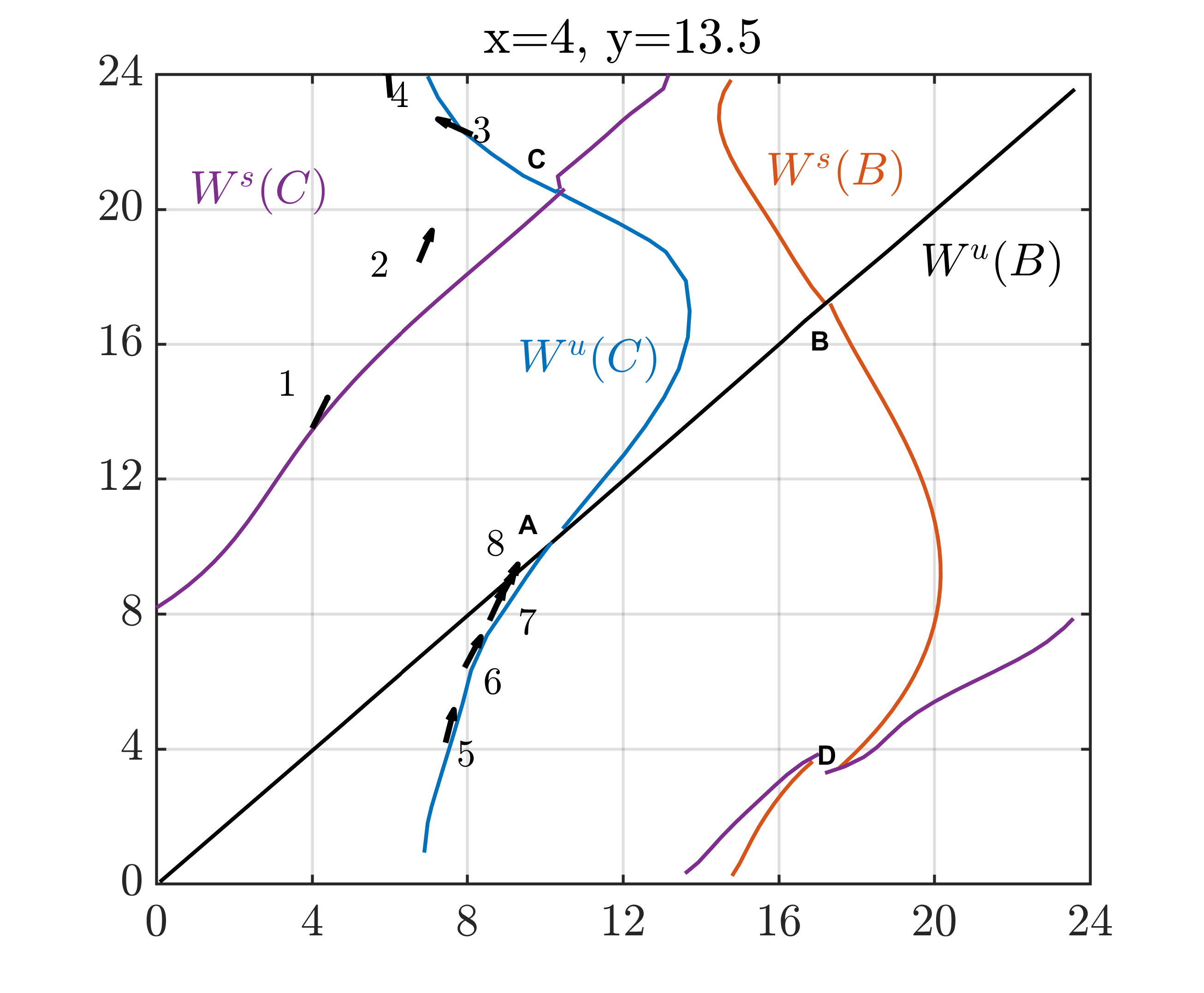}
			\includegraphics[width=.5\linewidth]{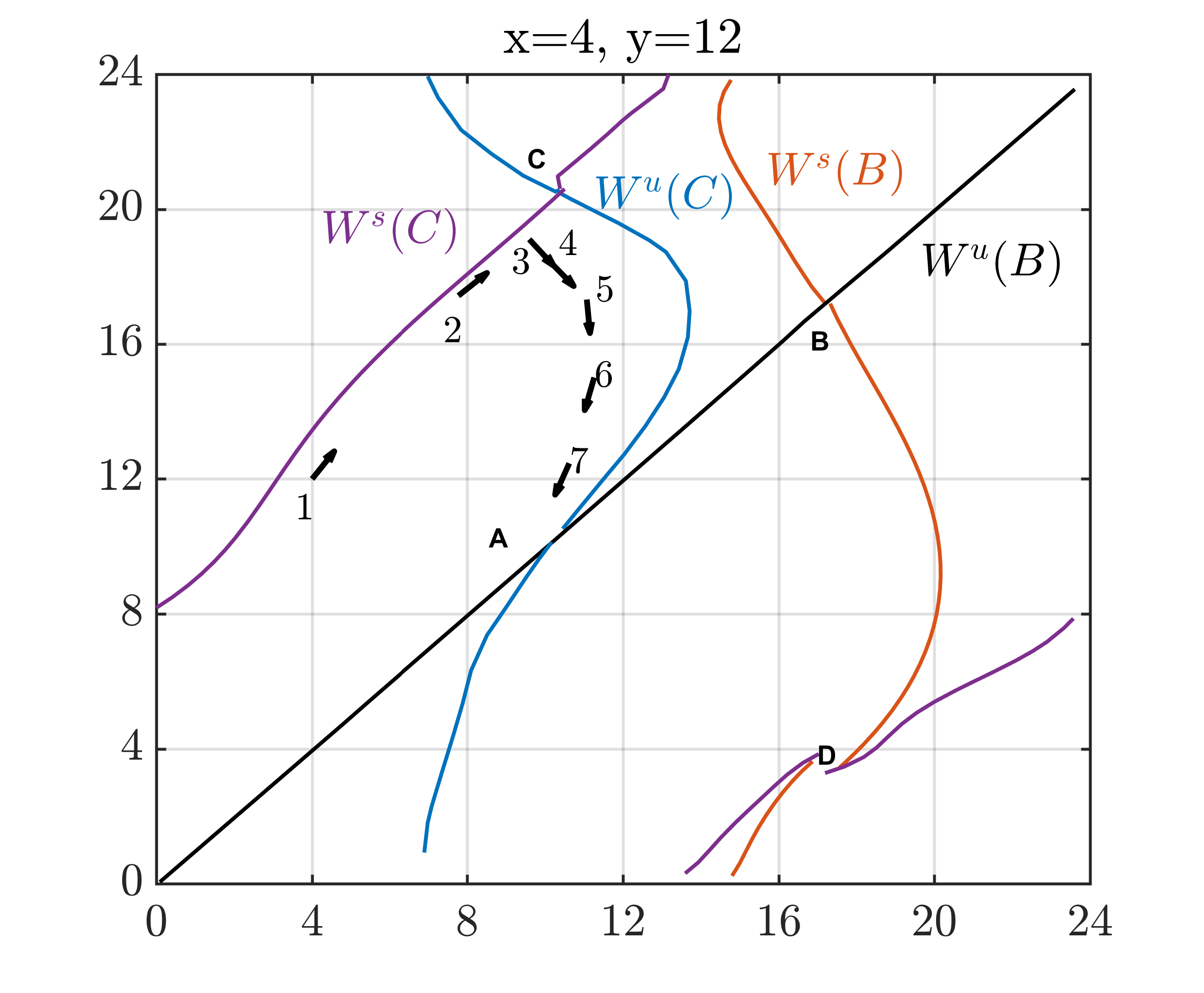}
			\caption{}
		\end{subfigure}
		\begin{subfigure}{.8\linewidth}
			\includegraphics[width=\linewidth]{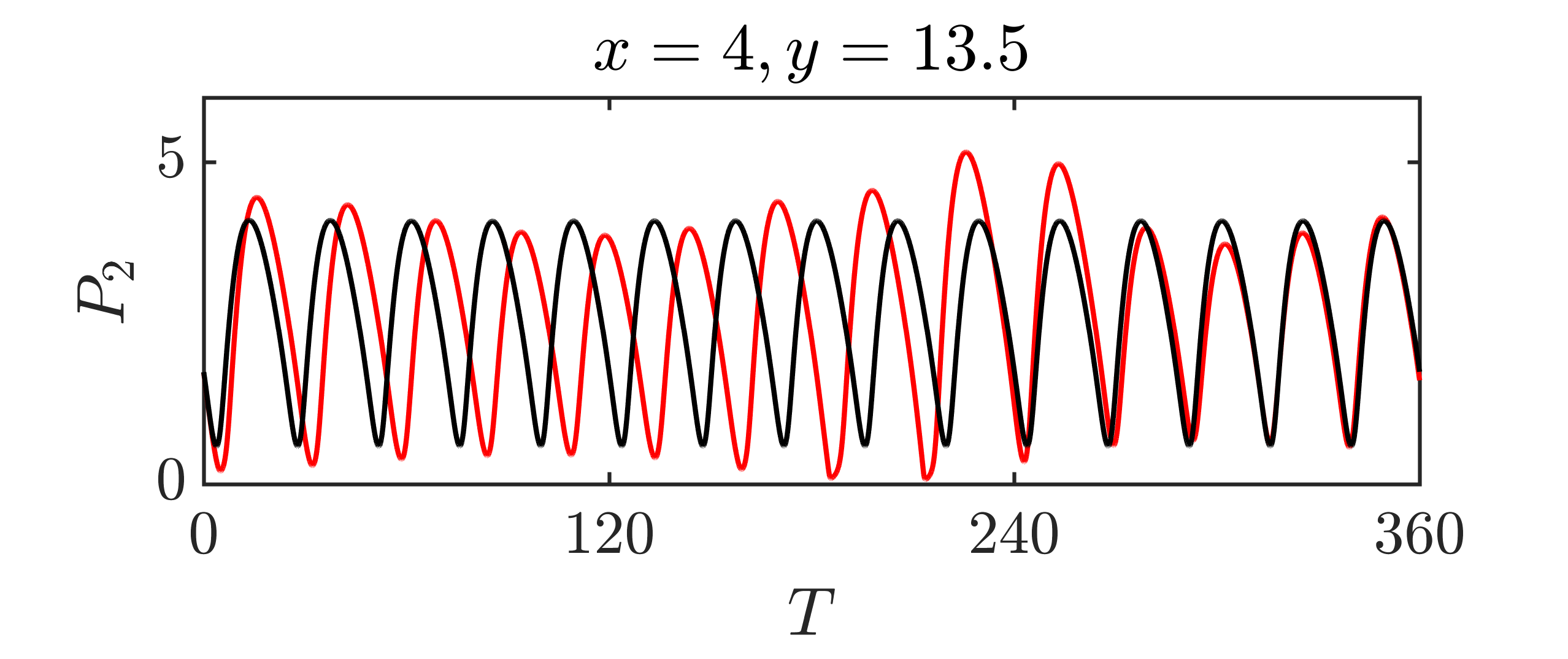}
			\includegraphics[width=\linewidth]{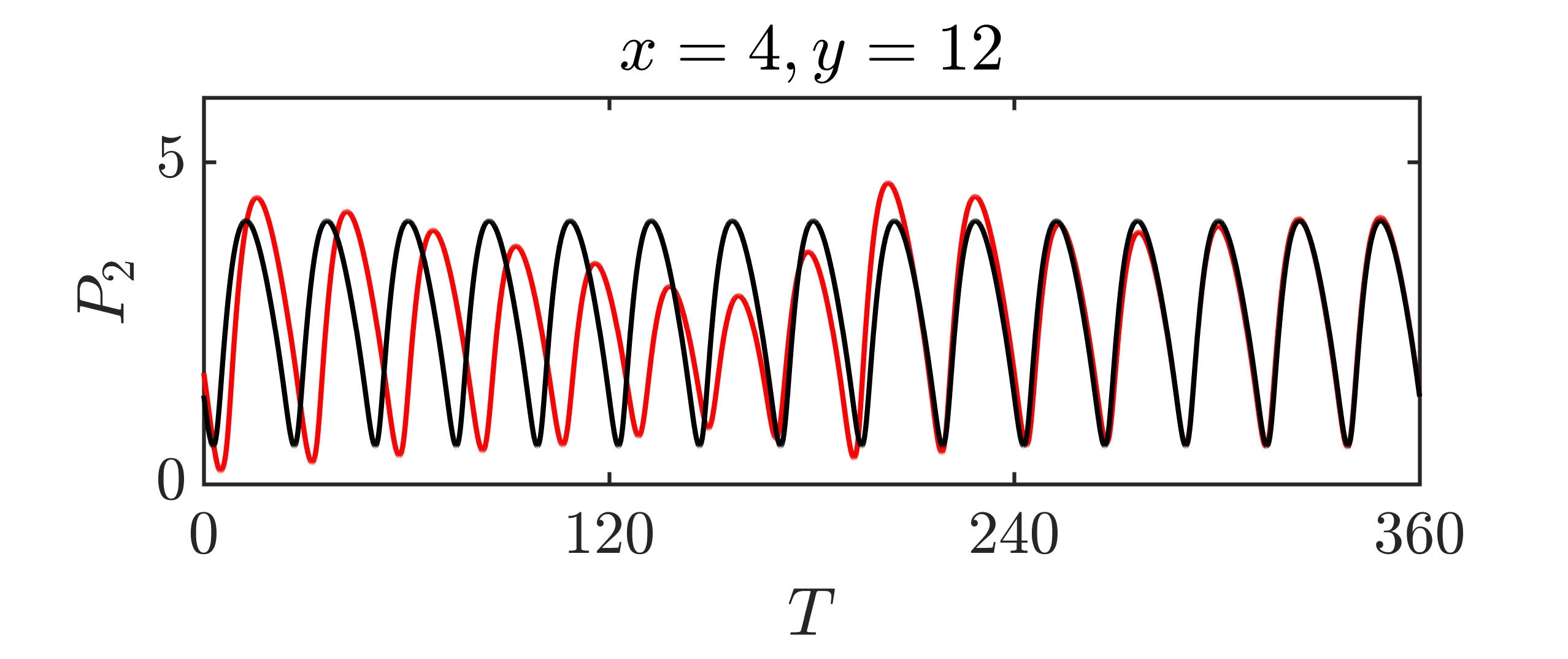}
			\caption{}
		\end{subfigure}
		\caption{Direction of entrainment depends sensitively on initial conditions. (a) The initial point (labeled 1) in the left panel lies above $W^s(C)$, while the similarly labeled point in the panel to the right lies below $W^s(C)$. Numbers indicate iterates. As shown, the direction of entrainment differs significantly. (b) Corresponding simulations agree with the iterates. Note the top panel shows that $O_2$ (red time course)  entrains through phase delay to the entrained solution (black time course); the lower panel shows $O_2$ entraining through phase delay-advance.}
		\label{fig:Iterates}
	\end{figure}

	The located manifolds are also helpful for understanding the direction of entrainment of 2-D maps. In the case of 1-D map, the direction of entrainment is essentially either phase advance or delay, and the longest entrainment times happen for initial conditions lying near the unstable fixed point. In the case of 2-D map, the direction of entrainment need no longer be monotonic. The manifolds associated with the saddle points B and C appear to behave like a separatrix, despite this being a map and not a flow. To classify the direction of entrainment in the 2-D map, we consider phase delays and advances in the $x$ and $y$ directions separately. For the $x$ direction, if the rotated angle from $x_n$ to $x_{n+1}$ is greater than $2\pi$, we call it phase delay, otherwise we call it phase advance. For the $y$ direction, we use the same definition as in the $O_1$-entrained map. To illustrate different directions of entrainment, we pick several initial conditions near the stable manifolds, then iterate the map. We also run simulations with the same initial conditions for comparison. For Fig. \ref{fig:Iterates}(a), in the left panel, we pick an initial point slightly above $W^s(C)$. It entrains to the stable solution by phase delay in the $y$ direction, and phase delay-advance-delay in the $x$ direction. In the right panel, the initial point is slightly below $W^s(C)$, but the entrainment is through phase delay-advance in $y$, and phase delay-advance in $x$. The corresponding simulations in Fig. \ref{fig:Iterates}(b) agree with the direction of entrainment found through the map and demonstrate the sensitivity to initial conditions.
	
	\paragraph{Parameter dependence of the map}
	
\begin{figure}[H]
		\centering
		\begin{subfigure}{.45\linewidth}
			\includegraphics[width=\linewidth]{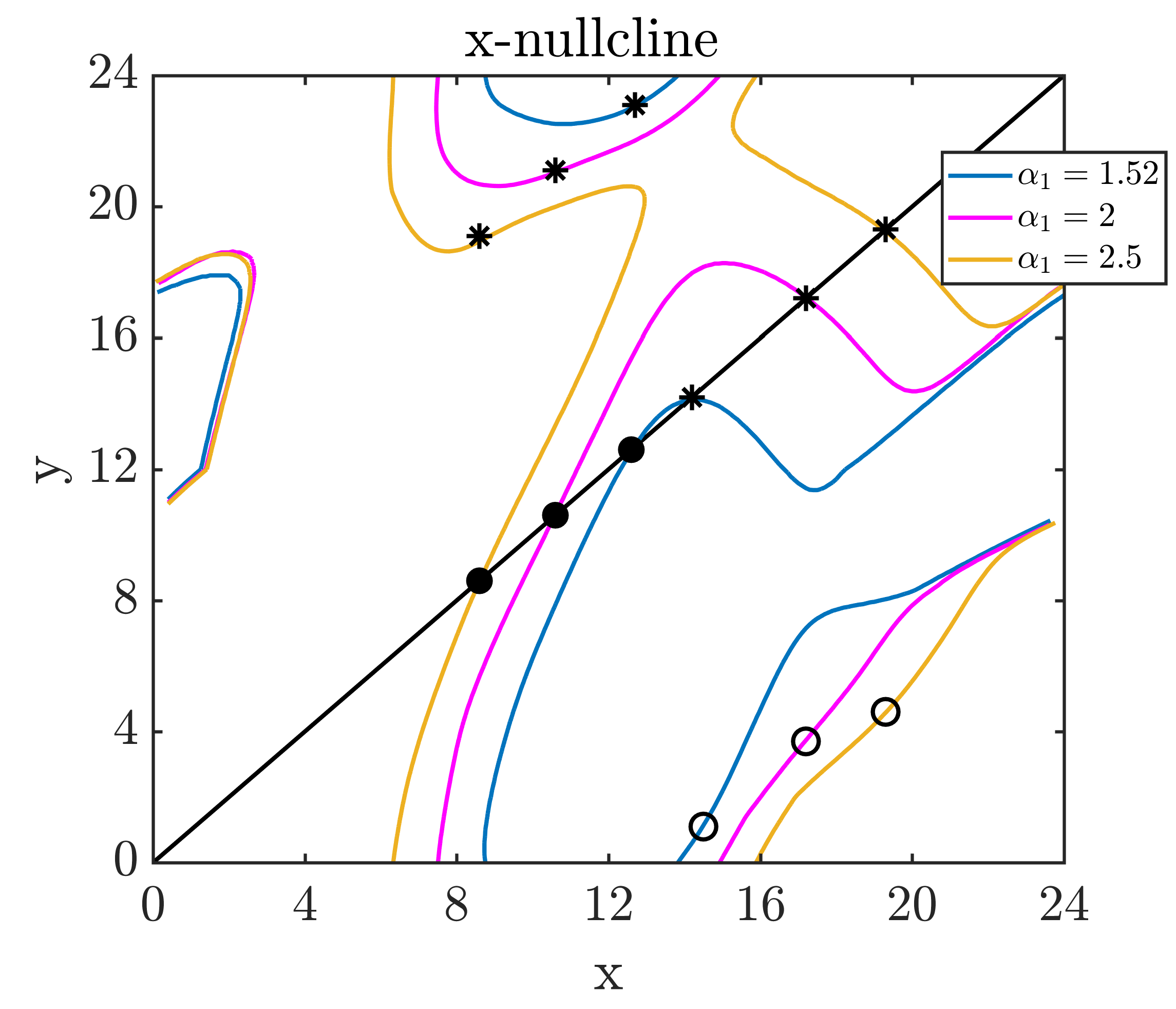}
			\caption{}
			\label{fig:xNullcline}
		\end{subfigure}
		\begin{subfigure}{.45\linewidth}
			\includegraphics[width=\linewidth]{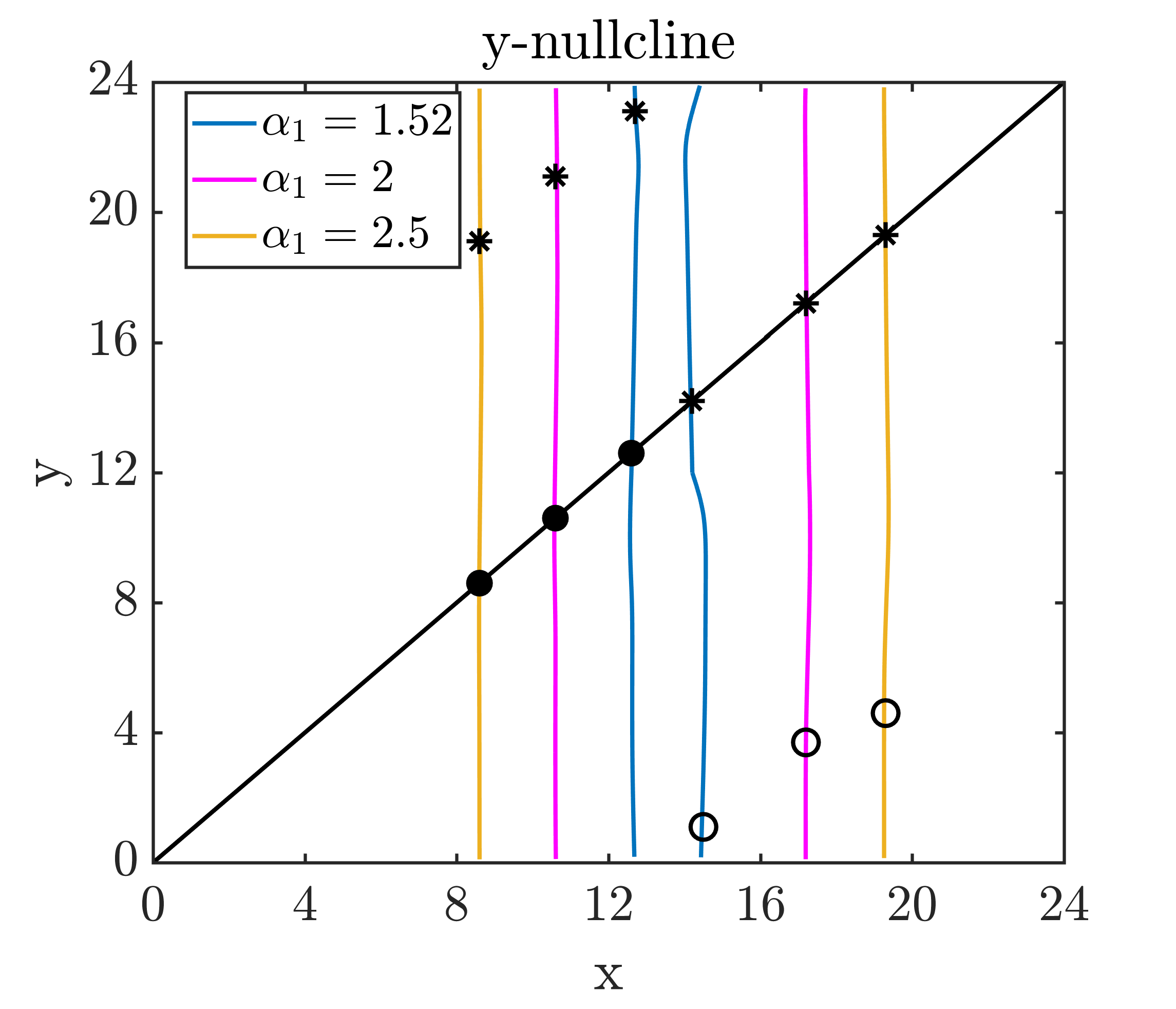}
			\caption{}
			\label{fig:yNullcline}
		\end{subfigure}
		\begin{subfigure}{.45\linewidth}
			\includegraphics[width=\linewidth]{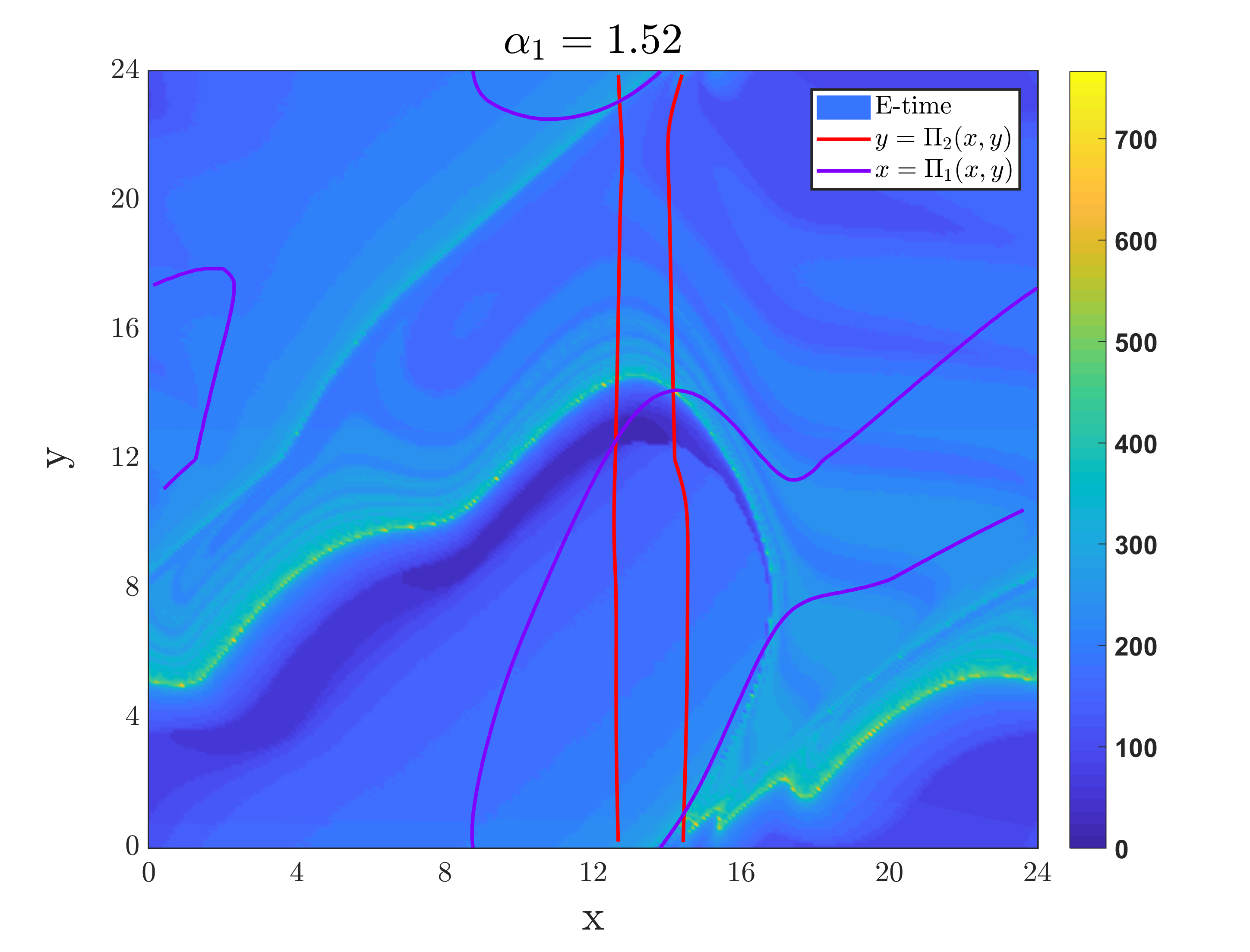}
			\caption{}
			\label{fig:Etime_alpha152}
		\end{subfigure}
		\begin{subfigure}{.45\linewidth}
			\includegraphics[width=\linewidth]{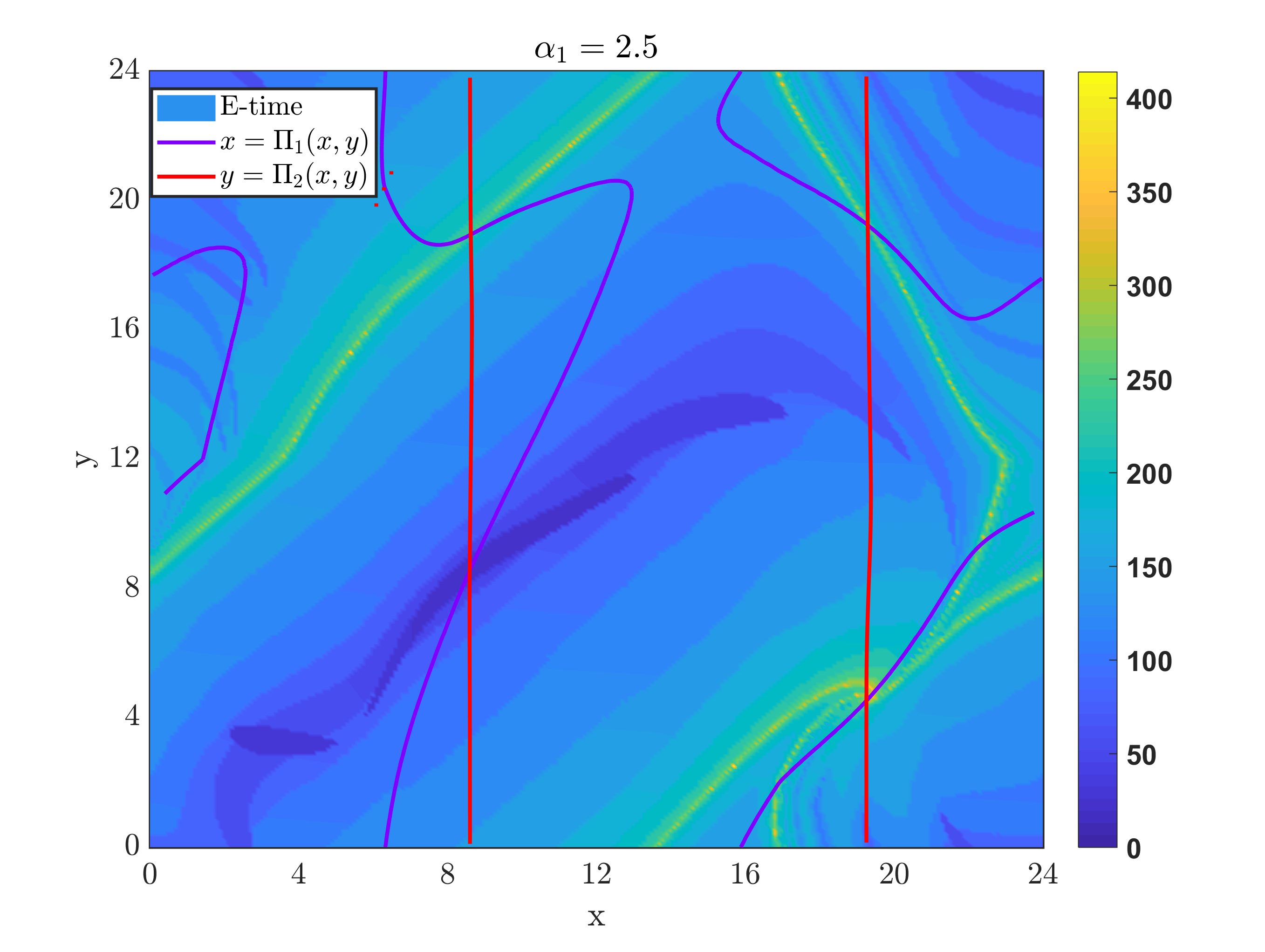}
			\caption{}
			\label{fig:Etimemap_alpha25}
		\end{subfigure}
		
		\caption{(a)-(b): The x and y nullclines under different $\alpha_1$ values. Solid circles denote unstable fixed points, open circles stable fixed points, and stars saddle points. (c)-(d): The heatmap of entrainment times for different values of $\alpha_1$. Note the difference in numeric value of the maximum value of the color scale.}
		\label{fig:6}
	\end{figure}

	In the section on the $O_1$-entrained map, we calculated the $O_1$-entrained map for four different values of $\alpha_1$, and found the system will lose entrainment if the coupling strength is too small. Now we calculate the 2-D map at different values of $\alpha_1$ to see how the fixed points and the entrainment time depend on $\alpha_1$. In Fig. \ref{fig:xNullcline} and \ref{fig:yNullcline}, we show the $x$ and $y$ nullclines for three different $\alpha_1$ values; the points with solid circle are the stable fixed points, the points with open circles are the unstable fixed points, and the starred points are saddle points. In Fig. \ref{fig:Etime_alpha152}, we show the heatmap of entrainment times for $\alpha_1=1.52$. In
	Fig. \ref{fig:Etimemap_alpha25}, we show the heatmap of entrainment times for $\alpha_1=2.5$. Note that $\alpha_1=2$ is our canonical case, and was presented before in Fig. \ref{fig:Etimemap_alpha_2}. Increasing $\alpha_1$, in general, decreases the entrainment time as can be observed from the color scale values (yellow max value $\approx 700$ for $\alpha_1=1.52$) verus $400$ for $\alpha_1=2.5$. In other words, stronger coupling between the central to peripheral oscillator speeds up entrainment.  
	
	\paragraph{The 2-D map for the semi-hierarchical case}
	For the strictly hierarchical model with only one feedforward connection from $O_1$ to $O_2$,  we have shown how to construct both the $O_1$-entrained map and the general 2-D entrainment map. Here we will show that the 2-D map can be derived for the model when $0<k_{L_2}<k_{L_1}$. In this case, $O_1$ is still dominant, allowing us to keep a semi-hierarchical structure.
	
\begin{figure}[H]
		\centering
		\begin{subfigure}{.45\linewidth}
			\includegraphics[width=\linewidth]{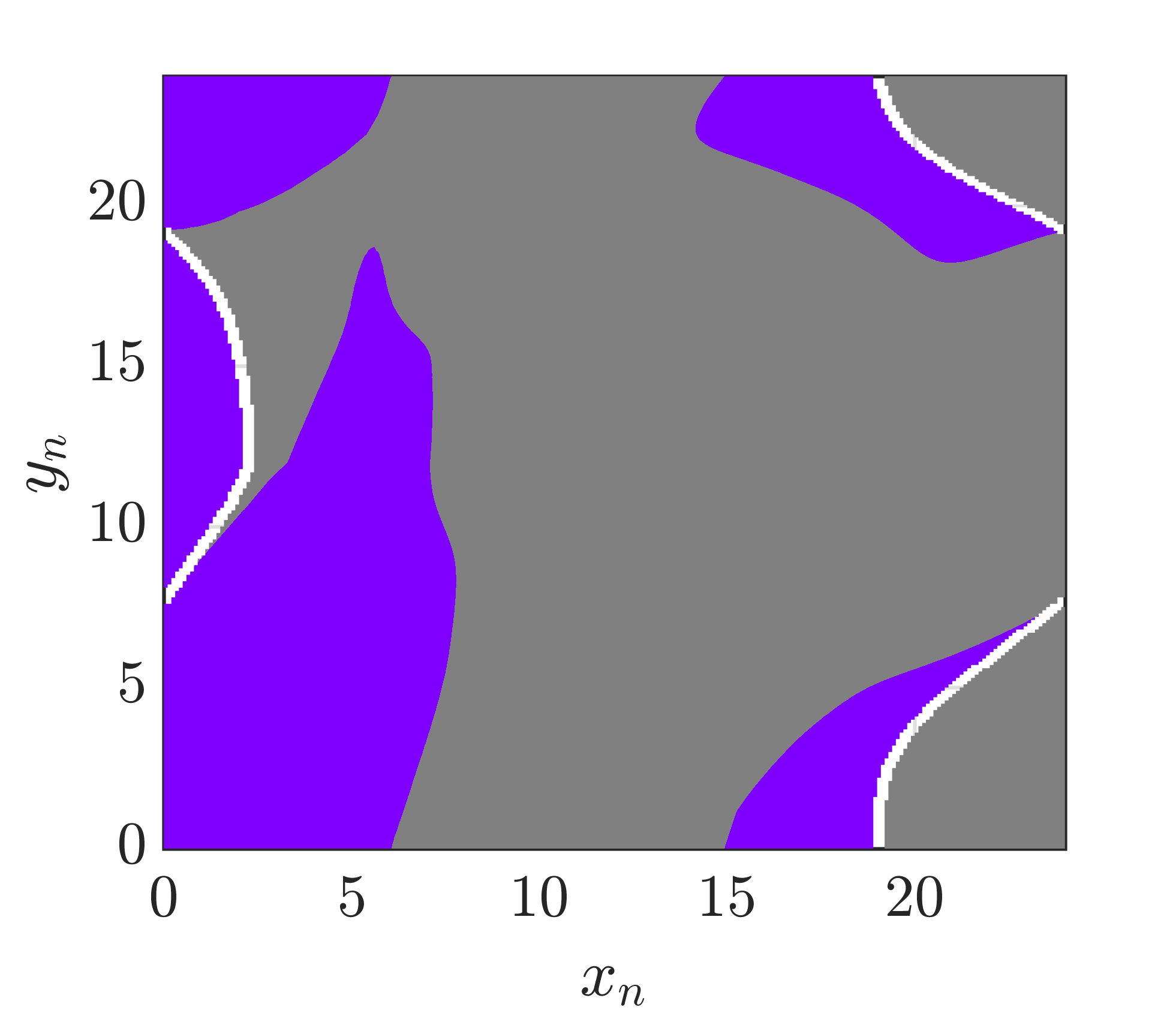}
			\caption{}
		\end{subfigure}
		\begin{subfigure}{.45\linewidth}
			\includegraphics[width=\linewidth]{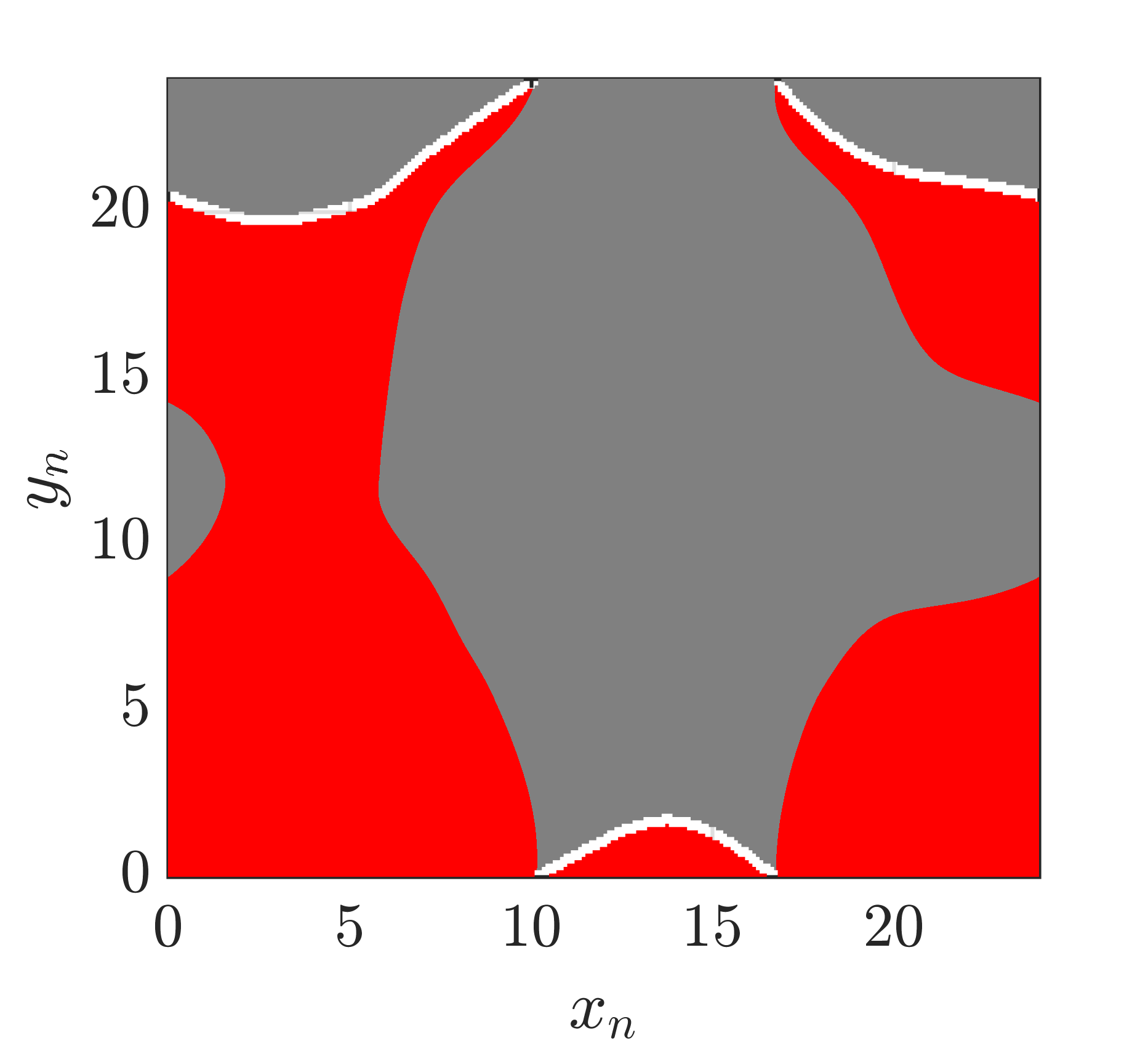}
			\caption{}
		\end{subfigure}
		\begin{subfigure}{.45\linewidth}
			\includegraphics[width=\linewidth]{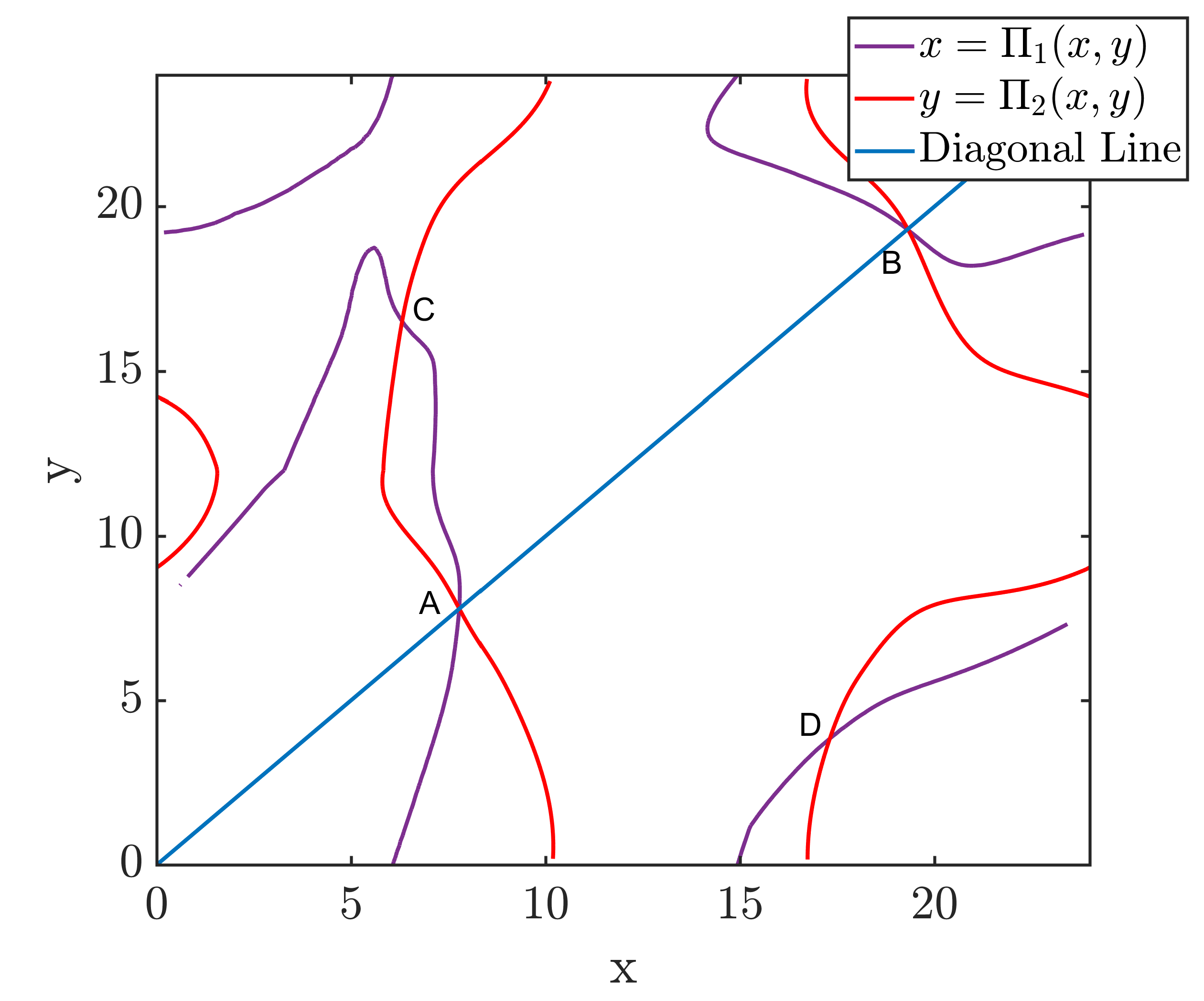}
			\caption{}
		\end{subfigure}
		\begin{subfigure}{.45\linewidth}
			\includegraphics[width=\linewidth]{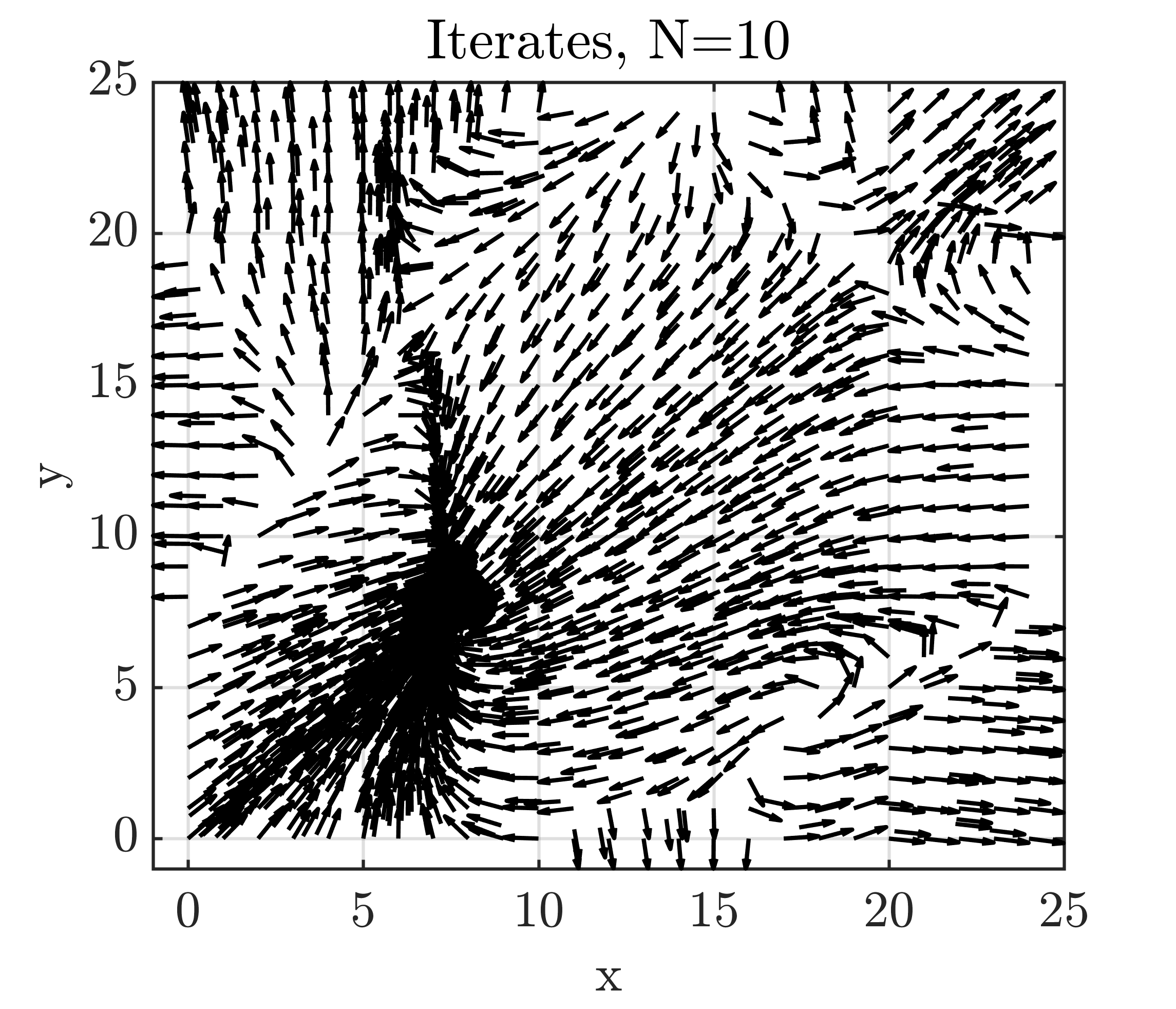}
			\caption{}
		\end{subfigure}
		
		\caption{2-D semi-hierarchical case. (a)-(b): The top view of $\Pi_1$ and $\Pi_2$ are presented; see Fig. \ref{fig:surf_pi1_top} and \ref{fig:surf_pi2_top} for an explanation of color coding. (c): we obtained 4 fixed points (A,B,C,D) with similar stability of the canonical model. (d): Ten iterates of each point. The vector field looks qualitatively similar to the strictly hierarchical case shown in Fig. \ref{fig:vf}a.}
		\label{fig:complete_model}
	\end{figure}

	We take $k_{L_2}=0.025$, and keep the values of other parameters the same, so that $O_1$ and $O_2$ both receive light forcing. We define the Poincar\'e section $\mathcal{P}: P_2=1.72, |M_2-0.1548|<\delta$.  We then obtained a 2-D map for this model. In Fig. \ref{fig:complete_model}(a) and (b), the top view of $\Pi_1$
	and $\Pi_2$ are presented. In Fig.~\ref{fig:complete_model}(c), we similarly obtained 4 fixed points (A,B,C,D) as in the strictly hierarchical case. Compare to the strictly hierarchical model, we found that the additional light forcing into $O_2$ accelerates the entrainment process, so that the time to return to $\mathcal{P}$ is decreased. Thus the whole surface shifts down, which causes A to move to the left of the diagonal, and B to move to the right of the diagonal. For points C and D, the limit cycle of $O_2$ is now determined by both $O_1$ and the light forcing, which changes the location of C and D. In Fig. \ref{fig:complete_model}(d), we calculated the first 10 iterates of each initial point. Comparing these results with the strictly hierarchical case, the stability of each fixed point remains unchanged, but their location has changed. Further, the entrainment time required for each initial condition is reduced because of the LD forcing into $O_2$.
	
	\section{Discussion}
	\label{sec:conclusion}
	
	Circadian oscillations exist from the sub-cellular level involving genes, proteins and mRNA up to whole body variations in core body temperature. These oscillations are typically entrained to the 24-hour light dark cycle. Additionally, food, exercise, exterior temperature and social interactions can also act as entraining agents in certain species \cite{mistlberger2005nonphotic}. In these cases, various pathways in each species exist which carry the entraining information to relevant parts of the circadian system. In this paper, we refer to the set of oscillators that first receive this input as central circadian oscillators. In turn, these central oscillators send signals about the time of day to other peripheral oscillators. When viewed in this manner, we obtain a hierarchical circadian system. For example, in the strictly hierarchical model (Fig.~\ref{fig:case1}), the central oscillator $O_1$ could represent the suprachiasmatic nucleus (SCN), the master pacemaker in the hypothalamus of mammals. The peripheral oscillator $O_2$ that does not receive light input could represent circadian clocks in organs such as the heart or kidney. Alternatively, $O_1$ could represent the part of the SCN that directly receives light input (the ventral core), and $O_2$ could then represent the part of the SCN that does not (the dorsal shell) \cite{honma2018mammalian}. For the semi-hierarchical model (Fig.~\ref{fig:case_nofeedback}), $O_1$ and $O_2$ could represent the central and peripheral clocks in {\it Drosophila}, since in flies the clock protein cryptochrome is a photoreceptor and thus even peripheral organs receive some direct light input \cite{chauhan2017central}.  The main goal of this paper has been to develop a low-dimensional method to study the basic properties of hierarchical systems such as the existence and stability of entrained solutions, together with how the phase and direction of entrainment of the constituent oscillators depends on important parameters. 
	
	In this work, we have developed a systematic method to assess whether and how hierarchical circadian systems entrain to external 24-hour light-dark forcing. With this approach, we were able to determine that over a large set of parameters, the 2-D map possesses four fixed points, each of which corresponds to a periodic orbit of the hierarchical circadian system. Only one of these fixed points is asymptotically stable. The other three fixed points are unstable. We showed how one of them, labeled D in Figs. \ref{fig:contours}, \ref{fig:vf}b and \ref{fig:complete_model}c, is a source from which iterates emerge, including the stable manifolds of the two saddle points B and C. These manifolds appear to act as separatrices in the $x$-$y$ domain of the map in the sense that, although they are for a map and not a flow, the manifolds separate the direction of convergence towards the stable fixed point A. Perhaps this is not so surprising as the saddle structure of the fixed points implies the existence of a saddle structure of the periodic orbits associated with points B and C. In the full five-dimensional phase space of the flow,  each of the corresponding one-dimensional stable and unstable manifolds from the map become three dimensional; the motion along the $O_1$ and $O_2$ limit cycles provide the additional two dimensions. This would be enough to form a separatrix in the five-dimensional phase space. The stable manifolds of the map have another effect. Any iterates that start close to either one of them lead to very long entrainment times. Consequences of this are discussed below.
	
	\paragraph{Related work}
	The mechanisms of communication between clock neurons is a topic of much ongoing research in the circadian field. The neuropeptide pigment-dispersing factor (PDF) is thought to act as the main synchronizing agent in the fly circadian neural network (Lin et al \cite{lin2004neuropeptide}). The analogue of PDF in the mammalian circadian system is vasoactive intestinal peptide (VIP), which plays a major role in synchronizing SCN neurons (Aton and Herzog \cite{mazuski2015circadian}). Although it is clear from studies with mutants that these neuropeptides provide important signals to synchronize circadian cells, the manner in which the signal interacts with the molecular clock is not well-understood (Dubrowy and Sehgal \cite{dubowy2017circadian}). Mathematical modeling can be used to explore the effect of different coupling mechanisms on clock network synchronization. In our model, we have assumed that production of the synchronizing factor is induced by activation of the clock gene in oscillator 1 ($M_1$), and that the effect of the synchronizing factor is to directly increase transcription of the clock gene is oscillator 2 ($M_2$). This type of coupling is similar to how Gonze et al \cite{GonzeHerzel05} modeled the the action of VIP in the mammalian clock network, however in the Gonze model they included a linear differential equation for the production and decay of the coupling agent. Thus, in their model the coupling agent is a delayed version of the clock gene activity. In the Roberts et al's model \cite{roberts2016functional} of the fly clock network, the coupling signal is also increased by clock gene activity. As in our model, the coupling signal then instantaneously increases the clock gene transcription rate in other oscillators. In addition, the Roberts model included a second type of coupling where the coupling signal depends on clock protein levels, rather than clock gene, and the effect of the coupling signal is to instantaneously reduce the clock gene transcription rate in other oscillators. Their simulations suggested that networks with both coupling types promoted synchrony and entrainment better than networks with either type of coupling alone. In a more detailed model of the fly clock network, Risau-Gusman and Gleiser \cite{risau2014mathematical} explored 21 different coupling mechanisms and found that synchronization of the network can only be achieved with a few of them. In future work, it would be interesting to use generalized entrainment maps to try to gain insight into why certain types of coupling promote synchrony and entrainment better than others.\\
	
	Several prior modeling studies on entrainment of circadian oscillators exist. Bordyugov et al \cite{bordyugov2015tuning} used the Kuramoto phase model and found, via Arnold tongue analysis, that the forcing strength and the oscillator amplitude both affect the entrainment speed. As noted in their work, a limitation of the method is that it only works for relatively weak coupling. An et al \cite{an2013neuropeptide} found that large doses of VIP (vasoactive intestinal polypeptide) reduce the synchrony in the SCN, which then reduces the amplitude of circadian rhythms in the SCN. In turn, they show that this leads to faster reentrainment of the oscillators in a jet lag scenario. Lee et al \cite{lee2017experimental} directly introduced a linear phase model to study the entrainment processes. They found that the period of the central and peripheral oscillators are not the only predictors of the entrained phase. The intensity of light forcing to the central oscillator and the strength of coupling from the central to the peripheral oscillator also play a role in determining the stable phase. Their results are consistent with what we found for the $O_1$-entrained map shown in Fig. \ref{fig:1DPEMap}. Roberts et al \cite{roberts2016functional} studied a population of coupled, modified, heterogeneous Goodwin oscillators under DD and single light pulse conditions.  Their model simulations of a semi-hierarchical system show that because of heterogeneity, a single light pulse can desynchronize and phase disperse the oscillators. This can lead to a change in the coupling strength between oscillators which in turn leads to a new periodic solution of different amplitude than before the light pulse. Although they didn't consider 24-hour LD forcing, Roberts et al suggest that this desynchrony can be an important component in assessing reentrainment of semi-hierarchical networks after jet lag. Our 2-D entrainment shows that this is indeed true. Namely, a shift in the light phasing that retains synchrony between $O_1$ and $O_2$ is equivalent to changing the initial $y$-value of our map, but keeping $x$ fixed. Whereas a shift of light phasing accompanied by a desynchronization is equivalent to changing both $x$ and $y$ from the stable fixed point. As our simulations show (Fig. \ref{fig:complete_model}d), the reentrainment process can be quite different in these two cases.

	There are two modeling papers of hierarchical systems that are quite relevant to our work. In Leise and Siegelman \cite{LeiseSiegelmann06}, the authors consider a multi-stage hierarchical system to assess properties of jet lag. They utilized a two dimensional circadian model due to olde Scheper et al \cite{Olde1999mathematical} to show that the direction of entrainment of peripheral oscillators need not follow that of the central oscillator. This is referred to as reentrainment by partition. To understand this idea more clearly, consider the concepts of orthodromic and antidromic reentrainment which are studied in the context of a time zone shift as in jet lag.  Orthodromic reentrainment is defined as the oscillator shifting in the same direction as the forcing signal (e.g. advancing in response to an advance of the light/dark cycle) and antidromic reentrainment is when the oscillator shifts in the opposite direction as the forcing signal (e.g. delaying in response to an advance of the light/dark cycle). The situation is more complicated for hierarchical systems where different parts of the system may shift in different directions. For example, Leise and Siegelman studied various scenarios involving 6 hour phase advances. They showed that a sudden change of this type leads to reentrainment though partition, but that if one spreads out the 6 hour advance over a few days, say 1.5 hours a day for 4 days, then the reentrainment process was orthodromic. In our paper that would be equivalent to starting with an initial $y$ value that is 6 less than the value of the stable fixed point at A (10.2-6=4.2), or systematically moving the $y$ value back by 1.5 after each iterate. For the parameter values of our paper, a 6 hour shift would  cause both oscillators to phase advance to the stable phase, which is similar to the orthodromic result in the Leise and Siegelman paper. However, a phase delay of the lights by say 10 hours or so would place the initial condition in the vicinity of the saddle fixed point C, leading to reentrainment through partition depending on the exact location relative to C. Thus our findings can be used to infer that the Leise-Siegelman multistage model also possesses unstable saddle fixed points whose properties govern the reentrainment process. A second more recent paper due to Kori et al \cite{kori2017accelerating} developed a hierarchical Kuramoto model to study the entrainment of circadian systems. They applied the model to predict the reentrainment time after two types of phase shifts, a single eight-hour shift versus a two-step shift with 4-hour shifts in each step. It turns out the latter requires fewer days to recover. In our paper, this can be related to the properties of stable manifolds of B or C. For example, in Fig.~\ref{fig:Etimemap_alpha_2}, for a single eight-hour shift near the fixed point A, the new point will stay close to $W^s(C)$, which makes the reentraiment time longer. For two successive four-hour shifts, the new point will be further from $W^s(C)$, which decreases the reentrainment time. This result generalizes findings from Diekman and Bose \cite{diekman2018reentrainment} and Kori et al. \cite{kori2017accelerating}.
	
	Regarding the numerical methods that we used to find stable and unstable manifolds,  we basically applied the search circle for stable manifolds \cite{england2004computing} and the growing method \cite{krauskopf1998growing} for unstable manifolds.  One difference between those methods and ours is the domain of the map, $\mathbb{R}^2$ versus a torus $\mathbb{T}^2$ in our case.  Instead of growing one curve, our manifold is cut off when it hits the boundary of the domain. We then restart the calculation at the equivalent periodic point of the domain, e.g. $x=24$ is reset to $x=0$.  Another difference is the terminating criteria for both growing and SC methods rely on calculating the arclength to a predetermined length.  However, in our map, the manifolds are generated from a certain point (the source D or the sink A), thus our algorithm terminates when those manifolds enter a neighborhood of the corresponding fixed points D and A.
	
	Recently Castej\'on and Guillamon derived a different 2-D entrainment map \cite{castejon2020phase}. This map applies to a single oscillator (not necessarily a circadian oscillator), subject to pulsed periodic input. The variables of their map are the phase and amplitude of the oscillator. They use phase-response curve type methods to show that their 2-D map is more accurate in tracking the phase-locking dynamics as compared to a 1-D map of simply phase. While they use the term 2-D entrainment map, it appears that their method applies to a class of problems that are different than the ones considered in this paper.
	
	\paragraph{Advantages and disadvantages of our method}
	The methods derived in this paper have the following advantages. Aside from allowing us to calculate entrainment times and directions as discussed above, the method provides a clear geometric description of why these results arise. Namely, the unstable manifolds of various fixed points organize the iterate structure of the dynamics. Our method does not specifically require the LD forcing to be weak in amplitude or short in duration. This is in contrast to methods that use phase response curves and thus require weak coupling or short duration perturbations \cite{brown2004phase,nakao2016phase}.
	
	Secondly, the dimension on which we perform analysis is significantly reduced from five to two dimensions. The classical Poincar\'e map can reduce the dimension of the original system by one. For example, Tsumoto et al. \cite{tsumoto2006bifurcations} construct a Poincar\'e map for 10-dimensional Leloup and Goldbeter model of the {\it Drosophila} molecular clock \cite{leloup1998model}, reducing the dimension to 9. The phase reduction techniques of Brown et al. \cite{brown2004phase} can reduce the dimension of limit cycle oscillators to 1-D, however this method is not accurate for strong coupling.
	
	There are some disadvantages of the map. First, the map only works to study local behavior near the stable limit cycle solutions. This is because we restrict the type of perturbations that we are considering to allow only for a shift of the LD cycle or a shift of the central oscillator along its own limit cycle. In particular, we don't know if there is an unstable or stable structure outside the basin of attraction of the stable entrained solution without additional analysis. Secondly, the phase angle method works well with two-dimensional systems. For higher dimensional systems, it would require additional assumptions. 
	
	\paragraph{Open questions and future directions}
	This work is based in part on analysis and in part on simulations. We have not proved that the correspondence of the findings of the 2-D map, e.g. existence and stability of fixed points, actually exist for the hierarchical system of ODEs. We would like to use a one-dimensional phase model, for example a Kuramoto model \cite{Kuramoto84} for each oscillator, to see if this proof can be made. 
	Alternatively, we believe this method of mapping should be applicable to other models, such as Goodwin \cite{Goodwin65}, Gonze \cite{GonzeHerzel05} or Forger, Jewett, Kronauer \cite{forger1999simpler} oscillators which are all higher dimensional. Verifying this, at the moment, would have to rely on checking agreement with simulations. The 2-D entrainment map should also be applicable to understand the interaction of circadian and sleep-wake rhythms to generalize the findings of Booth et al \cite{booth2017one}.
	
	A necessary condition of our method is the existence of limit cycle solutions of the forced system, so that we can map any point in the phase plane to a point on the limit cycle. Light input is not the only forcing signal that a circadian oscillator receives. For instance, exercise, the intake of meals and taking melatonin can also be considered as an external forcing. We would like to develop the entrainment map for multiple forcing signals.
	Another possible direction for future work involves generalization of model reduction techniques. Most reduction techniques are based on weak coupling, such as phase reduction \cite{brown2004phase}. We would like to develop a technique for a system with strong coupling. This part could potentially be done by deriving a Floquet normal form \cite{castelli2015parameterization} in phase and amplitude space. For the circadian oscillators that we studied in this paper, it remains open how to derive the Floquet normal form.

	\noindent
	{\bf Acknowledgements:} The work of Casey Diekman was supported, in part, by  NSF grant DMS-1555237 and U.S. Army Research Office grant W911NF-16-1-0584.
	
	\bibliographystyle{siamplain}
	
	\bibliography{CNT_reference}
	
	\clearpage
	
\end{document}